\documentclass[11pt]{amsart}

\usepackage{amsxtra}
\usepackage{amssymb}

\newcommand\n{\nu}

\newcommand\dd{{\mathfrak d}}

\newcommand\rr{{\mathfrak r}}
\newcommand\kk{{\mathfrak k}}
\newcommand\hh{{\mathfrak h}}

\newcommand\nn{{\mathfrak n}}
\newcommand\ggo{{\mathfrak g}}

\newcommand\ee{{\mathfrak e}}
\newcommand\aff{{\mathfrak{aff}}}

\newcommand\CC{\mathbb C}

\newcommand\RR{\mathbb R}

\newcommand\HH{\mathbb H}
\newcommand\DD{\mathbb D}

\newcommand\Aut{\operatorname{Aut}}
\newcommand\ad{\operatorname{ad}}

\theoremstyle{plain}

\newtheorem{thm}{Theorem}[section]
\newtheorem{lem}[thm]{Lemma}
\newtheorem{prop}[thm]{Proposition}
\newtheorem{cor}[thm]{Corollary}

\theoremstyle{definition}

\newtheorem{remark}[thm]{Remark}
\newtheorem{example}[thm]{Example}

\begin{document}
\title[Invariant pseudo K\"ahler metrics]{Invariant pseudo
K\"ahler metrics in dimension four}

\author{Gabriela Ovando}
\address{CIEM - Facultad de Matem\'atica, Astronom\'\i a y F\'\i sica,
Universidad Nacional de C\'or\-do\-ba, C\'or\-do\-ba~5000, Argentina}
\email{ovando@mate.uncor.edu}

\date{\today}

\begin{abstract} Four dimensional simply connected Lie groups admitting a
pseudo K\"ahler metric are determined. The corresponding Lie algebras are
modelized and the  compatible pairs $(J,\omega)$ are parametrized up to complex
isomorphism (where $J$ is a complex structure and 
$\omega$ is a symplectic structure). Such structure gives rise to a pseudo 
Riemannian metric $g$,  for which $J$ is a parallel. It is proved that most 
of these complex homogeneous spaces admit a  pseudo K\"ahler Einstein 
metric. Ricci flat and flat metrics are determined. In particular Ricci flat
unimodular  K\"ahler Lie algebras are flat in dimension four. Other 
algebraic and geometric features are treated. A general construction  of Ricci 
flat pseudo K\"ahler structures in higher dimension on some affine Lie algebras 
is given. Walker and hypersymplectic metrics on Lie algebras are compared. \end{abstract}

\thanks{{\it (2000) Mathematics Subject Classification}: 32Q15, 32Q20, 53C55, 32M10, 57S25, 22E25 }
\keywords{ pseudo K\"ahler metrics, dimension four, invariant metrics, K\"ahler Lie algebras}

\maketitle

\section{Introduction}
Simply connected Lie groups endowed with a left invariant pseudo Riemannian K\"ahler metric are in correspondence with K\"ahler Lie algebras. 
K\"ahler Lie algebras are Lie algebras $\ggo$ endowed with a pair $(J, \omega)$ consisting of a complex structure $J$ and  a compatible symplectic structure $\omega$. 
A K\"ahler structure on a Lie algebra determines a pseudo-Riemannian metric $g$ defined as
$$g(x,y) = \omega (Jx, y) \qquad x,y \in \ggo$$
not necessarily definite, and for which $J$ is parallel. The Lie algebra $(\ggo, J, g)$ is also known as a Pseudo-K\"ahler Lie algebra or indefinite K\"ahler Lie algebra.
 K\"ahler Lie algebras are special cases of symplectic Lie algebras and of pseudo metric Lie algebras and therefore tools of both fields can be used to their study. 

Lie algebras (resp. homogenous manifolds) admitting a definite K\"ahler metric 
were exhaustive studied by many authors. Indeed the
condition of the pseudo-metric to be definite imposes restrictions on the
structure of the Lie algebra  \cite{BG2} \cite{DN} \cite{DM} \cite{LM}.
In the nilpotent case the metric  associated to a pair $(J, \omega)$ 
cannot be definite positive \cite{BG1}. However this is not the case in  general for solvable Lie algebras.

In this paper we describe K\"ahler four dimensional Lie algebras.  
Since four dimensional symplectic Lie algebras must be solvable \cite{Ch}, our
results concern all possibilities in this dimension. Similar studies in the six dimensional 
nilpotent case were recently given in \cite{CFU2}. 

We prove that four dimensional completely solvable K\"ahler Lie algebras and
$\aff(\CC)$ are modelized by one of the  following short exact sequences of 
Lie algebras:
$$ \begin{array}{ll}
0 \longrightarrow \hh \longrightarrow \ggo \longrightarrow J\hh\longrightarrow 0
 \qquad & \mbox{orthogonal sum}\\ \\
0 \longrightarrow  \hh \longrightarrow \ggo \longrightarrow \kk \longrightarrow 0 \qquad & \hh \mbox{ and } \kk 
\mbox{ $J$-invariant subspaces }\end{array}$$
where in both cases $\hh$ is an $\omega$-lagrangian ideal on $\ggo$ and 
(hence abelian) and $J\hh$ and $\kk$ are $\omega$-isotropic subalgebras. The
first sequence splits and the second one does not necessarily splits. There are
also three kind of non completely solvable four dimensional Lie algebras 
admitting a K\"ahler structure. 
In all cases the  compatible pairs $(J,\omega)$ are parametrized  up to complex 
isomorphism. 

 The geometric study of these spaces continues by writing the corresponding
 pseudo Riemannian metric.   We give the explicit computations of  the corresponding Levi Civita 
connection, curvature  and Ricci tensors which could be used for further 
proposes. Making use of these information and the  models we find totally 
geodesic submanifolds. Moreover it is proved that the neutral metric on the 
Lie algebras satisfying the second short exact sequence is a  Walker metric on $\ggo$.

 We prove that in 8
  of the  11 families of K\"ahler Lie 
  algebras there exists an Einstein representative
  among the compatible pseudo K\"ahler metrics.  
  
  We also determine all Ricci flat metrics. On the one hand we show the
  equivalence in the unimodular case between Ricci flat and flat metrics in
   dimension four. On the other hand  we prove that in dimension four, aside 
   from the hypersymplectic Lie
 algebras \cite{Ad}, any Ricci flat metric is provided either 
 by $(\RR \times \ee(2),J)$, with $\ee(2)$ the Lie algebra of 
 the group of rigid motions of $\RR^2$ or by $(\aff(\CC),J_2)$, the real Lie 
 algebra underlying the Lie algebra of the affine motions group of $\CC$.
 Furthermore the Ricci flat pseudo metrics are deformations of flat pseudo
 K\"ahler metrics. Hence in dimension four, K\"ahler Lie algebras admitting Ricci
 flat pseudo K\"ahler metrics are in correspondence with K\"ahler Lie algebras
 with flat pseudo Riemannian metrics.
 

If we look at the Lie algebras admitting abelian complex structures we prove 
that a Lie algebra which admits  this kind of complex structure and is 
symplectic is also K\"ahler. Moreover if this is the case, $(\ggo,J)$ is 
K\"ahler if and only if $J$ is abelian. For instance the Lie algebra 
$\aff(\CC)$ has both abelian and non abelian complex structures; however only the 
abelian ones admit a compatible symplectic form. 


Finally we generalize our  results  constructing K\"ahler 
structures on affine Lie algebras,  $\aff(A)$ , where $A$ is a
commutative algebra. This kind of Lie algebras cover all cases of four dimensional Lie algebras having abelian
complex structures \cite{BD2}. We give examples in higher dimensions of Ricci 
flat pseudo Riemannian metrics by generalizing the K\"ahler structure of $(\aff(\CC),J_2)$ to affine Lie 
algebras $\aff(A)$ where $A$ is a commutative complex associative algebra.
It is also proved that a Walker K\"ahler metric on a Lie algebra $\ggo$  can be
hypersymplectic whenever
 some extra condition is satisfied. In particular a Walker metric compatible with 
the canonical complex structure of $\aff(\CC)$ 
is shown.

In a final section we
compute the obtained pseudo Riemannian metrics in global coordinates.

All Lie algebras are assumed to be real along this paper. 

\section{Preliminaries}

K\"ahler Lie algebras are endowed with a pair $(J, \omega)$ consisting of a complex structure $J$ and  a compatible symplectic structure $\omega$: $\omega(Jx,Jy)=\omega(x,y)$, namely a K\"ahler structure on $\ggo$.

Recall that an {\it almost complex}  structure on a Lie algebra $\ggo$ is an endomorphism $J:\ggo \to \ggo$ satisfying $J^2=-I$, where $I$ is the identity map. The almost complex structure $J$ is said to be integrable if $N_J\equiv 0$ where $N_J$ is the tensor given by
\begin{equation}\label{NJ}
N_J(x,y)=[Jx,Jy]-[x,y] - J[Jx,y] - J[x,Jy] \qquad\mbox{ for all }  x,y \in \ggo.  
\end{equation}
An integrable almost complex structure $J$ is called a {\it complex structure} on $\ggo$.

An equivalence relation is defined among Lie algebras endowed with complex structures. The Lie algebras with complex structures $(\ggo_1,J_1)$ and $(\ggo_2,J_2)$ are equivalent if there exists an isomorphism of Lie algebras $\alpha:\ggo_1 \to \ggo_2$ such that $J_2 \circ \alpha = \alpha \circ J_1$.

Examples of special classes of complex structures are the abelian ones and those that determine a complex Lie bracket on $\ggo$.

 A complex structure $J$ is said to be {\it abelian} if it satisfies $[JX, JY] = [X,Y]$ for all $X, Y \in \ggo$. A complex structure $J$  introduces on $\ggo$ a structure of complex Lie algebra if $J \circ ad_X = ad \circ J X$ for all $X \in \ggo$, and so $(\ggo, J)$ is a {\it complex Lie algebra}, and that means that the corresponding simply connected Lie group is also complex, that is, left and right multiplication by elements of the Lie group are holomorphic maps.

A {\it symplectic structure} on a 2n-dimensional Lie algebra $\ggo$ is a  closed 2-form $\omega \in \Lambda^2(\ggo^{\ast})$ such that $\omega$ has maximal rank, that is, $\omega^n\neq 0$. 
 Lie algebras (groups) admitting symplectic structures are called {\it symplectic } Lie algebras (resp. Lie  groups).

The existence problem of compatible pairs $(J, \omega)$ on a Lie algebra $\ggo$ is set up to complex isomorphism. In other words to search for K\"ahler structures on $\ggo$ it is sufficient to determine the compatibility condition between any symplectic structure and a representative of the class of complex structures.
  In fact, assume that there is a complex structure $J_1$ 
   for which there exists a  symplectic structure $\omega$ satisfying  
  $\omega (J_1 X , J_1 Y) = \omega (X,Y)$ for all $X, Y \in \ggo$ and assume that  $J_2$ is equivalent to $J_1$. Thus there exists an 
  automorphism $\sigma \in \Aut(\ggo)$ such that $J_2 = \sigma_{\ast}^{-^1} J_1 \sigma_{\ast}$. Then it holds
$$\omega (X , Y) = {\sigma^{\ast}}^{-1}\sigma^{\ast} \omega (X , Y) = {\sigma^{\ast}}^{-1}\omega (J_1 {\sigma_{\ast}}X,J_1\sigma_{\ast}Y)
=\omega (J_2 X,J_2Y).$$

K\"ahler Lie algebras belong to the class of symplectic Lie algebras. 
Special objects on a symplectic Lie algebra $(\ggo,\omega)$ are the isotropic 
and lagrangian subspaces. Recall that a subspace $W\subset \ggo$  is called $\omega$-isotropic if and only if 
 $\omega(W,W) =0$ and is said to be  $\omega$-lagrangian if it is $\omega$-isotropic and $\omega(W,y)=0$ 
 implies $y \in W$.

\begin{lem} \label{tota} Let $(\ggo, J, \omega)$ be a K\"ahler Lie algebra. Then if $\hh$ is a isotropic ideal, then:

\begin{itemize}
\item $\hh$ is abelian

\item J($\hh$) is a isotropic  subalgebra of $\ggo$.
\end{itemize}
Thus $\hh + J\hh$ is a subalgebra of $\ggo$ and the sum is not necessarily 
direct. However $\hh \cap J\hh$ is an ideal of $\hh + J\hh$ invariant by $J$.

\end{lem}
\begin{proof} Since $\hh$ is a isotropic ideal, the first assertion follows 
from the condition of $\omega$ of being closed.

The integrability condition of $J$ restricted to $\hh$, which was proved to be 
abelian, implies
$$[Jx,Jy] = J([Jx,y] + [x, Jy])$$ showing that $J\hh$ is a subalgebra of 
$\ggo$. The compatibility between $J$ and $\omega$ says that 
$\omega(Jx,Jy)=\omega(x,y)=0$ for $x,y\in \hh$, and so $J\hh$ is isotropic.
 Furthermore if $\hh$ is $\omega$-lagrangian, then $J\hh$ is
 $\omega$-lagrangian, and the second assertion is proved.
\end{proof}

A K\"ahler structure on a Lie algebra determines a pseudo-Riemannian metric $g$ defined as
\begin{equation}\label{metric}
g(x,y) = \omega (Jx, y) \qquad x,y \in \ggo
\end{equation}
for which $J$ is parallel with respect to the Levi Civita connection for $g$. Note that $g$ is not necessarily definite; the signature  is $(2k,2l)$ with 
$2(k+l)= \dim \ggo$.
 
Conversely if $(\ggo,J,g)$ is a Lie algebra endowed with a complex structure $J$ compatible with the pseudo metric $g$ then (\ref{metric}) defines a 2-form compatible with $J$ which is closed if and only if $J$ is parallel \cite{KN}. 
Hence the Lie algebra $(\ggo, J,g)$ is called a K\"ahler Lie algebra with pseudo-(Riemannian) K\"ahler metric $g$. 

Let $g$ be a pseudo Riemannian metric on $\ggo$. For a given  subspace $W$ of $\ggo$, the orthogonal subspace $W^{\perp}$ is
defined as usual by
$$W^{\perp} = \{ x \in \ggo / g(x,y) = 0,\, \text{ for all } y\in W\}.$$
The subspace $W$ is said to be isotropic if $W \subset W^{\perp}$ and is called totally 
isotropic if $W=W^{\perp}$.

 Lemma (\ref{tota}) can be rewritten in terms of the pseudo-Riemannian metric $g$.

\begin{lem} \label{iso} Let $(\ggo, J, g)$ be a K\"ahler Lie algebra. Assume that  an ideal $\hh \subset \ggo$ satisfies $J\hh\subset \hh^{\perp}$. Then

\begin{itemize}
\item $\hh$ is abelian and 

\item J$(\hh)$ is a subalgebra of $\ggo$ with $\hh \subset (J\hh)^{\perp} = J(\hh^{\perp}):=J \hh^{\perp}$.
\end{itemize}

Thus $\hh + J\hh$ is a subalgebra of $\ggo$ invariant by $J$ and the sum is not necessarily direct. However $\hh \cap J\hh$ is a $J$ invariant ideal of $\hh + J\hh$.
\end{lem}
\begin{proof} The subspace $\hh$ is $\omega$-isotropic if and only if $J \hh \subset \hh^{\perp}$. Hence $\hh$ is $\omega$-lagrangian if and only if $J\hh = \hh^{\perp}$. These remarks prove the assertions.
\end{proof}

 In \cite{DM} it is proved that if  $\ggo$ is a K\"ahler Lie algebra whose respective metric is positive definite then $\ggo$ is isomorphic to $\hh \rtimes J\hh$ when  $\ggo$ admits a ideal $\hh$ such that $\hh^{\perp}= J\hh$.

\subsection{On four dimensional solvable Lie algebras} It is known that a four 
dimensional symplectic Lie algebra must be solvable \cite{Ch}.  Let us recall the classification of four dimensional solvable real Lie algebras. For a proof see for instance \cite{ABDO}.
Notations used  along this paper are compatible with the following table. 

\begin{prop} \label{clases} Let $\ggo$ be a solvable four dimensional real Lie 
algebra. Then if $\ggo$ is not abelian, it is equivalent to one and only one 
of the Lie algebras listed below:
$$\begin{array}{ll}
{\rr\hh_3:} & {[e_1,e_2] = e_3 }  \\
{\rr\rr_3:\qquad } & {[e_1,e_2] =  e_2,\,[e_1,e_3] = e_2+ e_3 }  \\
{\rr\rr_{3,\lambda}:} & {[e_1,e_2] =e_2, [e_1,e_3]=\lambda e_3} \qquad {\lambda \in [-1,1]}\\
{\rr \rr'_{3,\gamma}:} & {[e_1, e_2] = \gamma e_2 -  e_3,[e_1, e_3] =  e_2 + \gamma e_3}  \qquad {\gamma \ge 0}\\
{\rr_2 \rr_2:} & {[e_1, e_2] = e_2,\,[e_3,e_4] = e_4}\\
{\rr_2':} & {[e_1,e_3] = e_3,\,[e_1,e_4] = e_4,\,[e_2,e_3]= e_4, \,[e_2,e_4]=-e_3} \\
{\nn_4:} & {[e_4,e_1] = e_2,\,[e_4, e_2] = e_3 } \\
{\rr_4:} & {[e_4,e_1]=e_1,\, [e_4,e_2] = e_1 + e_2, [e_4, e_3]=e_2 + e_3}\\
{\rr_{4,\mu}:} & {[e_4,e_1] = e_1,\, [e_4,e_2] = \mu e_2,\,[e_4,e_3] = e_2 + \mu e_3 } \qquad {\mu \in \RR} \\
{\rr_{4,\alpha, \beta}:}& {[e_4,e_1] = e_1,\, [e_4,e_2] = \alpha e_2,\, [e_4,e_3] = \beta e_3,}\,  \\
& \text{{\rm with }} \,-1 < \alpha \leq \beta \le 1 , \, \alpha \beta \ne 0,\,  \text{\rm{or }} \,-1 = \alpha \leq \beta\le 0 \\
{\rr'_{4,\gamma, \delta}:} & {[e_4,e_1] = e_1,\, [e_4,e_2] = \gamma e_2 - \delta e_3, \, [e_4,e_3] = \delta e_2 +\gamma e_3} \quad \gamma \in \RR, \delta > 0\\
{\dd_4:} & {[e_1,e_2]=e_3,\, [e_4,e_1] = e_1,\, [e_4,e_2] = -e_2}   \\
{\dd_{4,\lambda}:} & {[e_1,e_2]=e_3,\, [e_4,e_3] = e_3,\,[e_4, e_1]=\lambda  e_1,\,\,[e_4, e_2]=(1-\lambda) e_2} \quad  \lambda \ge \frac12\\
{\dd'_{4,\delta}:} & {[e_1,e_2]=e_3,\, [e_4, e_1]=\frac{\delta}2  e_1- e_2,[e_4,e_3] = \delta e_3,\,\,[e_4, e_2]=e_1+ \frac{\delta}2 e_2} \quad \delta \ge 0\\
{\hh_4} & {[e_1,e_2]=e_3,\, [e_4,e_3] = e_3,\,[e_4, e_1]=\frac{1}2 e_1,\,\,[e_4, e_2]= e_1 + \frac{1}2 e_2 }  \\ 
\end{array}
$$
\end{prop}

\begin{remark}
Observe that $\rr_2\rr_2$ is the Lie algebra  $\aff(\RR) \times \aff(\RR)$, where $\aff(\RR)$ is the Lie algebra of the Lie group  of affine motions of $\RR$, $\rr_2'$ is  the real Lie algebra underlying on the complex  Lie algebra $\aff(\CC)$, 
$\rr'_{3,0}$ is the trivial extension of $\ee(2)$, the Lie algebra of the  Lie group of rigid motions of $\RR^2$;  $\mathfrak r_{3,-1}$ is the Lie algebra $\mathfrak
e(1,1)$ of the group of rigid motions of the Minkowski $2$-space; $\rr\hh_3$ is the trivial extension of
the three-dimensional Heisenberg Lie algebra  denoted by $\hh_3$.
\end{remark}

A Lie algebra is called {\it unimodular} if tr($\ad_x$)=0 for all $x\in \ggo$, where tr denotes the trace of the map. 
 The unimodular four-dimensional solvable Lie algebras algebras  are: $\RR^4$, $\rr\hh_3$, $\rr \rr_{3,-1}$, $\rr\rr'_{3,0}$, $\nn_4$, $\rr_{4,-1/2}$, $\rr_{4,\mu,-1-\mu}$ $(-1<\mu\leq -1/2)$, $\rr'_{4,\mu,-\mu/2}$, $\dd _{4}$, $\mathfrak d_{4,0}'.$


Recall that a solvable Lie algebra is {\it completely solvable} when $\ad_x$ has real eigenvalues for all $x\in\ggo$.


Invariant complex structures in the four dimensional solvable real case were 
classified by J. Snow \cite{Sn} and G. Ovando \cite{O1}.   The following propositions show all Lie algebras of dimension four admitting special kinds of complex structures, making use of notations in (\ref{clases}).

\begin{prop}\label{cabel} If $\ggo$ is a  four dimensional Lie algebra admitting an abelian complex structure, then $\ggo$ is isomorphic to one of the following Lie algebras:  $\RR^4$,  $\RR \times \hh_{3}$, $\RR^2 \times  \aff(\RR)$ $\aff(\RR) \times \aff(\RR)$, $\aff(\CC)$, $\dd_{4,1}$.
\end{prop}
\begin{proof} If $\ggo$ is a four dimensional Lie algebra admitting abelian
complex structures then $\ggo$ must be solvable and its commutator has dimension
at most two (see \cite{BD2}). Let $\ggo$ be a four dimensional Lie algebra
satisfying these conditions. The first case is the abelian one which clearly
possesses an abelian complex structure. If $\dim \ggo'=1$ then $\ggo$ is
isomorphic either to $\RR\times \hh_3$ or to $\RR \times \aff(\RR)$, both
admitting abelian complex structures (see \cite{Sn} or \cite{BD1}). If the commutator 
is two dimensional then it must be abelian and therefore $\ggo$ must satisfy
 the following splitting short exact sequence of Lie algebras:
 $$ 0 \longrightarrow \RR^2 \longrightarrow \ggo \longrightarrow \hh \longrightarrow 0$$
with $\hh \simeq \aff(\RR)$ or $\RR^2$. The solvable four dimensional Lie algebras which satify these 
conditions are: $\RR^4$, $\rr \hh_3$, $\rr \rr_{3,\lambda}$, $\rr \rr_3$, 
$\rr \rr_{3,\lambda}'$, $\rr_2 \rr_2$, $\rr_2'$, $\dd_{4,1}$ (see for example 
\cite{ABDO}). The Lie algebras $\rr\rr_3$, $\rr_{3,\lambda}', \rr_{3,\lambda}$ 
$\lambda \neq 0$ do not admit abelian complex structures and  the other 
Lie algebras admit  such  kind of complex structures (see (\cite{Sn})).
\end{proof}

\begin{prop} Let $\ggo$ be a solvable four dimensional Lie algebra such that $(\ggo,J)$ is a complex Lie algebra, then $\ggo$ is either $\RR^4$ or   $\aff(\CC) = \rr_2'$.
\end{prop}
\begin{proof} Let $(\ggo,J)$ be Lie algebra with a complex structure $J$
satisfying $J[x, y]=[Jx, y]$ for all $x,y\in \ggo$. Then $J\ggo'
\subset \ggo'$
and hence $\dim \ggo'=$ 2 or 4. Assume now that $\ggo$ is solvable but not abelian and let
$x, Jx$ be a basis of $\ggo'$. Let $y, Jy$ not in $\ggo'$ such that $\{x,Jx,y,
Jy\}$ is a basis of $\ggo$. Then $[Jy,y]=0=[x,Jx]$ and the action of $y,Jy$
restricted to $\ggo'$  has the form
$$\ad_y=\left( \begin{matrix} a & -b \\ b & a \end{matrix} \right) \qquad 
\ad_{Jy}=\left( \begin{matrix} b & a \\ -a & b\end{matrix} \right)$$ where $a$
and $b$ are real numbers such that $a^2 +b^2\neq 0$. This implies that $\ggo\simeq \aff(\CC)$. In
fact taking $y'=\frac{1}{a^2+b^2}(a y+ b Jy)$ then $\{y', Jy', x, Jx\}$ is a
basis of $\ggo$ satisfying the Lie bracket relations of $\rr_2'$ in
(\ref{clases}).
 \end{proof}

 \section{Four dimensional K\"ahler Lie algebras}
In this section we determine all four dimensional K\"ahler Lie algebras 
and we  parametrize their compatible pairs $(J, \omega)$.

Most K\"ahler Lie algebras can be  found in a constructive way. In fact, 
according to \cite{O2} any symplectic Lie algebra $(g, J, \omega)$ which is 
either completely solvable or isomorphic to $\aff(\CC)$ admits a 
$\omega$-lagrangian ideal or equivalently in terms of the pseudo metric $g$ 
admits an ideal $\hh$ with $J\hh=\hh^{\perp}$.

In  four dimensional K\"ahler Lie algebras admitting such ideal $\hh$ there 
are two possibilities for $\hh \cap J\hh$:
  it is trivial or coincides with $\hh$. If it is trivial then $\ggo$ is 
  isomorphic to $\hh \rtimes J\hh$. Hence we have the following splitting 
  short exact sequence of Lie algebras
\begin{equation}\label{s1} 0 \longrightarrow \hh \longrightarrow \ggo \longrightarrow J\hh \longrightarrow 0.\end{equation}
If $\hh \cap J\hh$ is not trivial, then $J\hh = \hh$. So $\ggo$ can be 
decomposed as $\hh \oplus \kk$, where $\hh$ and $\kk$ are $J$-invariant totally 
isotropic subspaces  and one has the short exact sequence of Lie algebras, which
does not necessarily splits:
\begin{equation}\label{s2}0 \longrightarrow \hh\longrightarrow \ggo 
\longrightarrow \kk\longrightarrow 0.\end{equation}
In both cases $\hh$ is abelian (\ref{tota}) and therefore will be identified 
with $\RR^2$. 

These facts will help us to construct four dimensional K\"ahler Lie algebras. The results of the following propositions can be verified with Table (\ref{kahler}). 

\begin{prop} \label{m11} Let $(\ggo,J, g)$ be a four dimensional K\"ahler Lie 
algebra.   Assume that there exists an abelian ideal $\hh$ such that the 
following splitting exact sequence holds
$$ 0 \longrightarrow \hh  \longrightarrow \ggo \longrightarrow J\hh 
\longrightarrow 0.$$
where the sum is orthogonal. Then $\ggo$ is isomorphic to: 
$\RR^4$, $\RR \times \hh_3$,  $\RR^2 \times \aff(\RR)$, $\aff(\CC)$, 
$\aff(\RR) \times \aff(\RR)$, $\rr_{4,-1,-1}$, $\dd_{4,1}$,
 $\dd_{4,2}$, $\dd_{4,1/2}$.
\end{prop}
\begin{proof} Let $J$ be an almost complex structure on $\ggo$ compatible with 
the pseudo Riemannian metric $g$. The splitting short exact sequence (\ref{s1})
 is equivalent to one of the following short exact sequences of Lie algebras
\begin{equation}\label{s4}0 \longrightarrow \RR^2 \longrightarrow \ggo \longrightarrow \RR^2 \longrightarrow 0.\end{equation}
\begin{equation}\label{s5} 0 \longrightarrow \RR^2 \longrightarrow \ggo \longrightarrow \aff(\RR) \longrightarrow 0.\end{equation}
The pseudo metric $g$ restricted to $\hh$ defines a pseudo Riemannian metric 
on the Euclidean two dimensional ideal. On $\RR^2$ up to equivalence there exist
 two pseudo Riemannian metrics: the canonical one and the indefinite one of signature (1,1).

{\rm Case (\ref{s4}):} If $\ggo$ is a Lie algebra satisfying the first sequence 
(\ref{s4}) then the almost complex structure $J$  is integrable if and only if
it satisfies
\begin{equation}\label{e6}[Jx,y]=[Jy,x] \qquad\mbox{ for all }x, y \in \hh
\end{equation} 
 and $J$ is parallel with respect to the Levi Civita connection for $g$ if and 
 only if
\begin{equation}\label{e7}
g([Jx, z], y)= g([Jy,z], x)\qquad \mbox{ for all } x, y, z \in \hh
\end{equation}
 For the canonical metric with the conditions (\ref{e6}) and (\ref{e7}) one  
 gets the Lie algebras $\RR^4$, $\RR^2 \times \aff(\RR)$, 
 $\aff(\RR) \times \aff(\RR)$. For the neutral metric one gets the Lie algebras
  $\RR\times \hh_3$, $\RR^2 \times \aff(\RR)$,$\aff(\RR) \times \aff(\RR)$,
   $\aff(\CC)$ and $\dd_{4,1}$.

\medspace

{\rm Case (\ref{s5}):} If $\ggo$ is a Lie algebra satisfying (\ref{s5}) then 
the almost complex structure $J$ is integrable if and only if 
\begin{equation}\label{e8}e_2= [Je_1,e_2]-[Je_2,e_1]
\end{equation} 
where $span\{e_1,e_2\}=\hh \simeq \RR^2$, and $J$ is parallel with respect to the Levi Civita connection for $g$ if and only if
\begin{equation}\label{e9}
g(Je_2, Je_k)= g([Je_2, e_k], e_1)- g([Je_1,e_k], e_2)\qquad \mbox{ for } k=1, 2 \end{equation}
 For the canonical metric with the conditions (\ref{e8}) and (\ref{e9}) one 
 gets the Lie algebras $\dd_{4,1/2}$, $\dd_{4,2}$. For the neutral metric one 
 gets the Lie algebras  $\rr_{4,-1,-1}$, $\dd_{4,1/2}$, $\dd_{4,2}$.
\end{proof}

\begin{prop} \label{m22} Let $(\ggo,J, g)$ be a four dimensional K\"ahler Lie 
algebra.   Assume that there exists an abelian ideal $\hh$ such that the  
short exact sequence of Lie algebras (\ref{s2})
$$0 \longrightarrow \hh\longrightarrow \ggo \longrightarrow \kk
\longrightarrow 0$$
 holds, where $\hh$ and $\kk$ are $J$ invariant totally isotropic subspaces.  
 Then $\ggo$ is isomorphic to: $\RR \times \hh_3$, $\aff(\CC)$, $\rr_{4,-1,-1}$, $\dd_{4,1}$, $\dd_{4,2}$.
\end{prop}
\begin{proof} At the algebraic level, the  sequence (\ref{s2}) takes the form 
(\ref{s4}) or (\ref{s5}), where $\hh =span\{e_1,Je_1\} \simeq \RR^2$ and 
$\kk=span\{e_2, Je_2\} \simeq \RR^2$ in (\ref{s4}) or $\aff(\RR)$ in (\ref{s5}).
 Let $J$ be a complex structure on $\ggo$ and let $\omega$ be a 2-form 
 compatible with $J$. Then  $\omega$  being  closed is equivalent to:
$$\omega([e_2,Je_2], x) + \omega([x, e_2], Je_2) + \omega([Je_2,x], e_2)=0.$$
If (\ref{s2}) splits then for the case (\ref{s4}) one gets the Lie algebra 
$\aff(\CC)$, and in the case (\ref{s5}) one gets the Lie algebras 
$\rr_{4,-1,-1}$, $\dd_{4,1}$ and $\dd_{4,2}$. If (\ref{s2}) does not splits 
then one gets $\RR \times \hh_3$. 
\end{proof}

Notice that according to the four dimensional classifications of complex structures
\cite{Sn} \cite{O1} and symplectic structures \cite{O2} the non completely
solvable Lie algebras which could admit compatible pairs $(J, \omega)$ are 
$\RR \times \ee(2)$, $\rr'_{4,0,\delta}$, $\delta \neq 0$, and $\dd'_{4,\delta}$ with $\delta \neq 0$. These Lie algebras admit K\"ahler structures (see Table (\ref{kahler})) and moreover the Lie algebras $\RR \times \ee(2)$ and $\rr'_{4,0,\delta}$ satisfy the following splitting short exact sequence of Lie algebras
$$0 \longrightarrow \hh=J\hh \longrightarrow \ggo \longrightarrow \hh^{\perp} \longrightarrow 0$$
where $\hh$ is an abelian ideal but not a $\omega$-lagrangian ideal of $(\ggo, \omega)$.

Let $\ggo$ be a Lie algebra admitting a complex structure $J$ and let us denote by $\mathcal{S}_c(\ggo,J)$ the set of all symplectic forms 
$\omega$ that are compatible with $J$. Our goal now is to parametrize the elements of $\mathcal{S}_c(\ggo,J)$ where $\ggo$ is a four dimensional Lie algebra. In the previous paragraphs we found the Lie algebras $\ggo$ for which $\mathcal{S}_c(\ggo,J)\neq \emptyset$ for some complex structure $J$.

Denoting   $\{e^i\}$ be the dual basis on $\ggo^{\ast}$ of the basis $\{e_i\}$ on $\ggo$ (as in (\ref{clases})), we adopt the abbreviation
$e^{ijk\hdots}$ for $e^i \wedge e^j \wedge e^k \wedge \hdots$.

\begin{prop} \label{kahler} Let $\ggo$  be a K\"ahler Lie algebra, then $\ggo$ is isomorphic to one of the following Lie algebras endowed with complex and compatible symplectic structures listed as follows:

\medspace
\begin{center}
\begin{tabular}{|lll|}\hline
{$\begin{array}{c} \vspace{-.27cm}\\ \ggo \\  \vspace{-.27cm}
\end{array} $} & {\rm Complex structure} & {\rm Compatible symplectic 2-forms} \\ \hline
$ \rr \hh_3:$ & $J e_1 = e_2,\, J e_3 = e_4$ & ${\begin{array}{l} \vspace{-.27cm}\\ a_{13+24}(e^{13} + e^{24})+ a_{14-23} (e^{14}-e^{23})+\\ +a_{12}e^{12},\,
{a_{13}^2+a_{14}^2\ne 0} \\  \vspace{-.27cm}
\end{array} }$\\ \hline
{$\rr \rr_{3,0}:$} &  $ J e_1 = e_2,\, J e_3 = e_4$ & ${ \begin{array}{c} \vspace{-.27cm}\\a_{12} e^{12} + a_{34}e^{34},\, {a_{12}a_{34} \ne 0} \\  \vspace{-.27cm}
\end{array} }$\\ \hline
{$\rr \rr'_{3,0}:$} & {$J e_1 = e_4,\, J e_2 = e_3$}& {$
\begin{array}{c} \vspace{-.27cm}\\a_{14}e^{14} + a_{23}e^{23},\, a_{14}a_{23}\ne 0 \\ \vspace{-.27cm}
\end{array} $} \\ \hline
{$\rr_2 \rr_2:$} & {$J e_1 = e_2,\, J e_3 = e_4$}& {$ \begin{array}{c} \vspace{-.27cm}\\a_{12}e^{12} + a_{34}e^{34},\, a_{12}a_{34}\ne 0 \\ \vspace{-.27cm}
\end{array} $} \\ \hline
{$\rr'_2:$} & {$J_1 e_1 = e_3,\, J_1 e_2 = e_4$} & {$
\begin{array}{c} \vspace{-.27cm}\\
a_{13-24}(e^{13} - e^{24})+a_{14+23}(e^{14} + e^{23}),\\ a_{13-24}^2 +a_{14+23}^2 \ne 0
\end{array}$}   \\
& {$J_2 e_1 = -e_2,\, J_2 e_3 = e_4$} & {$\begin{array}{l} \\  a_{13-24}(e^{13} - e^{24})+a_{14+23}(e^{14} + e^{23})+\\ + a_{12} e^{12},\,
a_{13-24}^2 +a_{14+23}^2 \ne 0 \\ \vspace{-.27cm}
\end{array} $} \\ \hline
{$\rr_{4,-1,-1}:$} & {$Je_4 =  e_1,\, J e_2 = e_3$} &  {$
\begin{array}{l} \vspace{-.27cm}\\
a_{12+34}(e^{12}+e^{34}) + a_{13-24} (e^{13}-e^{24})+ \\
+a_{14} e^{14},\, a_{12+34}^2+a_{13-24}^2\ne 0 \\ \vspace{-.27cm}
\end{array}$} \\ \hline
{$\rr'_{4,0, \delta}:$} & {$
\begin{array}{l} \vspace{-.27cm}\\
J_1 e_4 = e_1,\, J_1 e_2 = e_3,\\
J_2 e_4 = e_1,\, J_2 e_2 = -e_3 \\ \vspace{-.27cm} \end{array} $} &
{$a_{14}e^{14} + a_{23}e^{23},\,
a_{14}a_{23}\ne 0 $} \\ \hline\end{tabular}

\begin{tabular}{|lll|} \hline
{$\dd_{4,1}:$} & {$J e_1 = e_4,\, J e_2 = e_3$} &{ $\begin{array}{l} \vspace{-.27cm}\\a_{12-34}(e^{12} - e^{34})+e_{14} e^{14},\, a_{12-34}\ne 0\\ \vspace{-.27cm} \end{array}$} \\ \hline
{$\dd_{4,2}:$} & {$\begin{array}{l} \vspace{-.27cm}\\
J_1e_4 = - e_2,\, J_1 e_1 = e_3 \\ \vspace{-.27cm} \\
J_2e_4 =  -2e_1,\, J_2 e_2 = e_3\\
\vspace{-.27cm} \end{array}$} & {$\begin{array}{l} \vspace{-.27cm}\\
a_{14+23}(e^{14} + e^{23})+ a_{24} e^2 \wedge e^4,\, a_{14+23}\ne 0 \\ \vspace{-.27cm} \\
a_{14} e^{14}+a_{23} e^{23},\, a_{14}a_{23}\ne 0\\
\vspace{-.27cm} \end{array}$}\\ \hline
{$\dd_{4,1/2}:$} & {$\begin{array}{l} \vspace{-.27cm}\\
J_1e_4 =  e_3,\, J_1 e_1 = e_2 \\
J_2e_4 =  e_3,\, J_2 e_1 = -e_2 \\
\vspace{-.27cm} \end{array}$} & {$\begin{array}{l} \vspace{-.27cm}\\
a_{12-34} (e^{12}-e^{34}),\, a_{12-34}\ne 0 \\
\vspace{-.27cm} \end{array}$} \\\hline
{$\dd'_{4,\delta}:$} & {$\begin{array}{l} \vspace{-.27cm}\\J_1e_4 =  e_3,\, J_1 e_1 = e_2,\\
J_2e_4 =  -e_3,\, J_2 e_1 = e_2,\\
J_3e_4 =  -e_3,\, J_3 e_1 = -e_2,\\
J_4e_4 =  e_3,\, J_4 e_1 = e_2,\\ \vspace{-.27cm} \end{array}
$}& { $\begin{array}{l}
a_{12-\delta 34} (e^{12}-\delta e^{34}),\, a_{12-34}\ne 0 \end{array}$} \\ \hline
\end{tabular}
\end{center}

\begin{center}
{\rm Table} \ref{kahler}
\end{center}

\end{prop}
\begin{proof} The complete proof follows a case by case study. Making  use of the classifications of complex structures we found in \cite{Sn} and \cite{O1}, then  for a fixed complex structure $J$ in a given Lie algebra $\ggo$ we verify the compatibility condition with the symplectic forms given in \cite{O2}.

We shall give the details in the case $\rr'_2$, the Lie algebra corresponding to $\aff(\CC)$. The other cases should be handled in a similar way. As we can see on the classification of Snow \cite{Sn} the complex structures on $\rr_2'$ are given by: $J_1 e_1 = e_3$, $J_1e_2 = e_4$;  and for the other type of complex structures, denoting $a_1\in \CC$ by $a_1 = \mu + i \nu$, with $\nu \ne 0$; we have   $J_{\mu,\nu} e_1 =\frac{\mu}{\nu} e_1 +(\frac{\nu^2 +\mu ^2}{\nu})e_2 $, $J_{\mu, \nu} e_3 = e_4$. On the other hand any sympletic structure has the form: $\omega = a_{12}(e^1 \wedge e^2) + a_{13-24}(e^1 \wedge e^3 - e^2 \wedge e^4)+a_{14+23}(e^1 \wedge e^4 + e^2 \wedge e^3)$, with $ a_{14+23}^2 + a_{13-24}^2 \ne 0$. Assuming that there exists a K\"ahler structure it  holds $\omega (JX,JY) = \omega (X,Y)$ for all $X,Y \in \ggo$ and this condition produces  equations on the coefficients of $\omega$ which should be verified in each case.

So for $J_1$ we need to compute only the following:
$$\omega (e_1,e_2) = a_{12} = \omega (e_3, e_4)$$
and $$\omega(e_1,e_4) = a_{14+23} = \omega (e_3, -e_2)$$
Thus these equalities impose the condition $a_{12}=0$. And so any K\"ahler structure corresponding to $J_1$ has the form  $\omega = a_{13-24}(e^1 \wedge e^3 - e^2 \wedge e^4)+a_{14+23}(e^1 \wedge e^4 + e^2 \wedge e^3)$ with $a_{13-24}^2+a_{14+23}^2\ne 0$.

For the second case corresponding to $J_{\mu,\nu}$, by computing $\omega(e_2,e_4)$, $\omega(e_1,e_3)$, we get respectively:
$$
\begin{array}{lrcl}
{\text{i)}} & (1 + \frac1{\nu}) a_{13-24} & = & \frac{\mu}{\nu} a_{14+23} \\
&&&\\
{\text{ii)}} & (1 + \frac{\mu^2 + \nu^2}{\nu}) a_{13-24} & = & -\frac{\mu}{\nu} a_{14+23} \\
\end{array}
$$
By comparing i) and ii)  we get:
$$
\begin{array}{lrcl}
 &(1 + \frac1{\nu})  a_{13-24} & = & -(1 + \frac{\mu^2 + \nu^2}{\nu}) a_{13-24}  \\
\end{array}
$$
and  this equality  implies either  iii)$ a_{13-24}=0$ or iv) $1 + \frac1{\nu}+ 1 + \frac{\mu^2 + \nu^2}{\nu} =0$. As $a_{13-24}\ne 0$ (since in this case we would also get $a_{14+23}=0$ and this would be  a contradiction) then it must hold iv), that is  $1+ \mu^2 + \nu^2 =-2\nu$ and that implies $ \mu^2 =-2\nu- 1- \nu^2= -(\nu +1)^2$  and  that is possible only if $\mu =0$ and $\nu = -1$. For this complex structure $J$, given by $Je_1=-e_2$ $ J e_3 = e_4$,   it is not difficult to prove  that for any symplectic structure $\omega$ it always holds $\omega(JX,JY) $ = $ \omega(X,Y)$, that is, any symplectic structure on $\ggo$ is compatible with $J$. In this way we have completed the proof of the assertion.
\end{proof}

In the following we shall simplify the notation: parameters with four subindices will be denoted only with two subindices, hence for instance $a_{14+23} \to a_{14}$. By the computations of the pseudo K\"ahler metrics the parameters satisfy those conditions of previos Table (\ref{kahler}).

\begin{remark} The dimension of $\mathcal{S}_c(\ggo,J)$ is 1, 2 or 3 in all of 
the cases. When $\dim \mathcal{S}_c(\ggo,J)= 1$, then $\mathcal{S}_c(\ggo,J)$ 
can be parametrized by $\RR^{\ast}$, when it is two, then by  
$\RR \times \RR^{\ast}$, $\RR^{\ast} \times \RR^{\ast}$  or $\RR^2-\{0\}$ and if
 it equals three by $\RR \times (\RR^2-\{0\})$.
\end{remark}

\begin{remark} \label{affc} The complex structure which endowes the Lie algebra $\rr_2'$ with
 a complex Lie bracket is given by $Je_1= e_2$ and $J e_3 = e_4$, which is not 
 equivalent with the listed  in the previous table. This complex structure does 
 not admit a compatible symplectic structure. In fact, assume that $\Omega$ is a
 2-form compatible with $J$, then $\Omega= \alpha e^{12} + \beta (e^{13}+ e^{24}) +
 \gamma (e^{14}- e^{23})$ . Hence $d\Omega = 0$ if and only if $\beta =0 =
 \gamma$. Thus there is no symplectic  structure compatible with $J$. 
 In \cite{CFU1} it is proved that any closed 2-form is always degenerate when 
 it is compatible with a  complex structure $J$ which gives $\ggo$ a structure
  of  complex Lie algebra.
\end{remark}

\begin{remark} Among the four dimensional Lie algebras we find many examples of Lie algebras, such that the set of complex structures  $\mathcal C$ and the set of symplectic structures $\mathcal S$ are both  nonempty and however there is no compatible pair $(J,\omega)$. This situation occurs for instance on the Lie algebras $\hh_4$ or the family $\dd_{4,\lambda}$ for $\lambda\neq 1/2, 1, 2$ (Compare results in \cite{O1} \cite{O2} and \cite{Sn}).
\end{remark}

Reading the previous list of Proposition \ref{kahler}
by looking at the structure of the Lie algebras we get the following Corollaries.
\begin{cor} Let $\ggo$ be a K\"ahler four dimensional  Lie algebra. If $\ggo$ is unimodular then it is isomorphic either to $\RR \times \hh_3$ or $\RR \times \ee(2)$.

If $\ggo$ is not unimodular then either:

i)  $\dim \ggo' = 1$ and it is isomorphic to $\RR^2\times \aff(\RR)$,

ii) $\dim \ggo' =2$ and $\ggo$ is a non trivial extension of $\ee(1,1)$, $\aff(\RR)\times\aff(\RR)$, or an extension of  $\ee(2)$ or

iii) $\ggo' \simeq \RR^3$ and $\ggo \simeq \rr_{4,-1,-1}$ or $\rr_{4,0,\delta}'$ or

iv) $\ggo' \simeq \hh_3$ and the action of $e_4\notin \ggo'$ diagonalizes with set of eigenvalues one of the following ones $\{1,1,0\}$, $\{1,2,-1\}$, $\{1,\frac12,\frac12\}$, $\{1,\frac12 + i \delta,\frac12-i\delta\}$, with $\delta >0$.
\end{cor}
\begin{proof} According to \cite{ABDO}, if $\dim \ggo' = 1$ then $\ggo$ is a
 trivial extension of $\hh_3$ or $\aff(\RR)$; the non trivial extension of $\ee(1,1)$ is $\rr_2\rr_2$ and the extensions of $\ee(2)$ are isomorphic either to $\aff(\CC)$ or $\RR \times \ee(2)$. The rest of the proof follows by looking at the adjoint actions on any K\"ahler Lie algebra with three dimensional commutator.
\end{proof}

\begin{cor} Let $\ggo$ be a nilpotent (non abelian) four dimensional K\"ahler 
Lie algebra, then it is isomorphic to $\RR\times \hh_3$ and any complex 
structure is abelian.
\end{cor}
\begin{proof} Among the four dimensional Lie algebras the non abelian nilpotent ones are $\RR \times \hh_3$ and $\nn_4$. Only $\RR\times \hh_3$  admits a compatible pair $(J, \omega)$ and in fact the previous table parametrizes elements of $\mathcal{S}_c(\RR \times \hh_3,J)$ for a fixed complex structure $J$.
\end{proof}
\begin{remark} $\RR \times \hh_3$ is the Lie algebra underlying the Kodaira Thurston nilmanifold \cite{Th} for which actually any complex structure $J$  admits a compatible symplectic form $\omega$.
\end{remark}

\begin{cor} Let $\ggo$ be a four dimensional Lie algebra for which any complex
structure gives rise to a K\"ahler structure on $\ggo$. Then $\ggo$ is
isomorphic either to $\RR \times \hh_3$, $\RR^2\times \aff(\RR)$, $\RR \times
\ee(2)$, $\rr_{4,-1,-1}$, $\rr_{4,0,\delta}'$, $\dd_{4,1}$ $\dd_{4,2}$ 
\end{cor}

\begin{cor} Let $\ggo$ be a four dimensional Lie algebra admitting abelian 
complex structures. Then  $(\ggo, J)$ is K\"ahler if and only if $\ggo$ is 
symplectic and $J$ is abelian.
\end{cor}
\begin{proof} According to (\ref{cabel}) and the results of (\cite{Sn}), the
four dimensional Lie algebras which are K\"ahler and admit abelian complex 
structures are $\RR \times \hh_3$, $\RR^2\times \aff(\RR)$, 
$\aff(\RR) \times \aff(\RR)$, $\aff(\CC)$ and $\dd_{4,1}$. Among these Lie 
algebras only $\aff(\CC)$ admits complex structures which are not abelian. 
On $\aff(\CC)$ there is a curve of non equivalent complex structures. Among 
the points of this curve there is one which belongs to the abelian class. The
class represented by this point and one class more corresponding to an
abelian structure admit a compatible symplectic structure and the complex
structure which are not abelian do not admit  a compatible symplectic 
structure.
\end{proof}

\

 In dimension four a  pseudo Riemannian K\"ahler
metric must be definite or neutral. Notice that the dimension of the set of 
pseudo Riemannian K\"ahler metrics on each K\"ahler Lie algebra $(\ggo, J, g)$ 
coincides with the dimension of $\mathcal{S}_c(\ggo,J)$.

We use the following notation to describe the pseudo Riemannian metrics. If 
$\{e_i\}$ is the basis of Proposition (\ref{clases}) then $\{e^i\}$ is its 
dual basis on $\ggo^{\ast}$ and symmetric two tensors are of the form $e^i\cdot e^j$ where $\cdot$ denotes the symmetric product of 1-forms.
We denote by $z_i$ the coordinates of $z\in\ggo$ with respect to the 
basis $\{e_i\}$.

\begin{cor}
Let $(\ggo,J)$ be a non abelian four dimensional K\"ahler Lie algebra with  
complex structure $J$ admitting only definite  K\"ahler metrics then $(\ggo,J)$ 
is isomorphic either to the Lie algebra $(\dd_{4,1/2},J_1)$, or to 
$(\dd_{4,\delta}', J_1, J_3)$.

The K\"ahler Lie algebras $(\RR \times \hh_3,J)$, $(\aff(\CC), J_1, J_2)$ , $(\rr_{4,-1,-1},J)$ and $(\dd_{4,1},J)$ and $\dd_{4,\delta}'$admit only neutral pseudo Riemannian metrics.
\end{cor}
\begin{proof} In the case of the completely solvable K\"ahler Lie algebras or 
$\aff(\CC)$ the assertions follow from the proof of Propositions (\ref{m11}) 
and (\ref{m22}). In fact these K\"ahler Lie algebras can be constructed in 
terms of splitting exact sequences of Lie algebras, verifying some extra 
conditions.   We need to study the assertions in the cases $\RR \times \ee(2)$,
 $\rr'_{4,0,\delta}$ and $\dd'_{4,\delta}$ with $\delta \neq 0$. Looking at the
 pseudo K\"ahler metrics on $\RR \times \ee(2)$,  $\rr'_{4,0,\delta}$ (see 
 Propositions (\ref{unif}) and (\ref{t4})) it is possible to verify  that  both cases admit  
 definite and neutral metrics. In the case of $\dd'_{4,\delta}$ the 
 complex structures $J_1$ and $J_3$ admit only definite compatible pseudo 
 metrics and the complex structures $J_2$ and $J_4$ admit only neutral 
 compatible pseudo metrics. 
\end{proof}

The following propositions offer an alternative model for four dimensional K\"ahler Lie algebras since  the existence of  a lagrangian ideal is a strong condition. The next constructions are based on the existence of an abelian ideal which does not need to be lagrangian.

\begin{prop} \label{m2} The following K\"ahler four dimensional Lie algebras: $(\RR^2 \times \aff(\RR), J)$, $(\RR \times \ee(2), J)$,  $(\aff(\RR) \times  \aff(\RR), J)$, $(\rr_{4,0,\delta}', J_1, J_2)$ endowed with a pseudo K\"ahler metric, satisfy  the following splitting short exact sequence of Lie algebras:
$$ 0 \longrightarrow \hh=J\hh \longrightarrow \ggo \longrightarrow \hh^{\perp} \longrightarrow 0 $$
where the sum is orthogonal.\end{prop}
\begin{proof} For the Lie algebras of the proposition, with a given pseudo K\"ahler metric, we exhibit a abelian ideal  satisfying  $J\hh=\hh$:
$$\begin{array}{lll}
\RR^2 \times \aff(\RR), J & g = a_{12}(e^1 \cdot e^1 +e^2 \cdot e^2)+a_{34} (e^3 \cdot e^3 + e^4 \cdot e^4) & \hh= spann\{e_1,e_2\}\\
\RR \times \ee(2), J & g = a_{14}(e^1 \cdot e^1 +e^4 \cdot e^4)+a_{23} (e^2 \cdot e^2 + e^3 \cdot e^3) & \hh= spann\{e_2,e_3\}\\
\aff(\RR)^2, J & g = a_{12}(e^1 \cdot e^1+ e^2 \cdot e^2)+a_{34} (e^3 \cdot e^3+ e^4 \cdot e^4) & \hh= spann\{e_1,e_2\} \\
\rr_{4,0,\delta}', J_1, J_2 & g = a_{14}(e^1 \cdot e^1+ e^4 \cdot e^4)+a_{23} (e^2 \cdot e^2+ e^3 \cdot e^3) & \hh= spann\{e_2,e_3\} \\
\end{array}
$$
\end{proof}

\begin{prop} \label{m3} The following K\"ahler four dimensional Lie algebras: 
$(\RR \times \hh_3, J)$, $(\rr_{4,-1,-1}, J)$, $(\dd_{4,2}, J_1)$,  endowed 
with a pseudo K\"ahler metric, satisfy  the following splitting exact sequence 
of Lie algebras:
$$ 0 \longrightarrow \hh=\hh^{\perp} \longrightarrow \ggo \longrightarrow J\hh \longrightarrow 0.$$
\end{prop}
\begin{proof} For the Lie algebras of the proposition, with a fixed pseudo K\"ahler metric $g$, we exhibit an  ideal  satisfying $\hh=\hh^{\perp}$ and $\hh \cap J\hh =0$:
$$\begin{array}{lll}
\RR\times \hh_3, J & g = e^1 \cdot e^3 - e^2 \cdot e^4 & \hh=spann\{ e_2,e_3\} \\
\rr_{4,-1,-1}, J & g= a_{13}(e^1 \cdot e^2-e^3\cdot e^4)  & \hh= spann\{e_1,e_3\} \\
\dd_{4,2}, J_1 & g = a_{14}( e^1 \cdot e^2 +  e^3 \cdot e^4)  & \hh= spann\{e_2,e_3\}\\
\end{array}
$$
\end{proof}
\begin{remark} The Lie algebras of Prop. (\ref{m3}) are those admitting an hypersymplectic structure \cite{Ad}.
\end{remark}

 \medspace

 \section{On the geometry of left invariant Pseudo K\"ahler metrics in  four dimensional Lie algebras}

 In this section we study  the geometry of the Lie group $G$ whose
  Lie algebra $\ggo$ is endowed with a K\"ahler structure.  
Because of the left invariant property all results in this sections are 
presented at the level of the Lie algebra. 
  We make use of the models (\ref{m11}) and (\ref{m22}) to find totally 
  geodesic submanifolds. We find Ricci flat and Einstein K\"ahler metrics. 
   In the definite case 
  Ricci flat metrics are flat \cite{A-K}. In the non definite case this is not 
  true in general. However in dimension four if $\ggo$ is unimodular and 
  the K\"ahler metric is Ricci flat, then it is flat.

\medspace

Let $\nabla$ be the Levi Civita connection corresponding to the pseudo Riemannian metric $g$. This is determined by the Koszul formula
$$ 2g(\nabla_x y, z) = g([x,y],z) -g([y,z],x) + g([z,x],y)$$
It is known that the completeness of the left invariant connection $\nabla$ on $G$ can be studied by considering the corresponding connection on the Lie algebra $\ggo$. Indeed the connection $\nabla$ on $G$ will be (geodesically) complete if and only if the differential equation on $\ggo$
$$
\dot x(t)=-\nabla_{x(t)} x(t)$$
admits solutions $x(t)\subset \ggo$ defined for all $t\in \RR$ (see for instance \cite{Gu}).

A submanifold $N$ on a Riemannian manifold $(M,g)$ is totally geodesic if $\nabla_x y \in TN$ for $x,y\in TN$. At the level of the Lie algebra we have totally geodesic subspaces, subalgebras, etc. which are in correspondence with totally geodesic submanifolds, subgroups, etc on the corresponding Lie group $G$ with left invariant pseudo  metric $g$.

\begin{prop} \label{l1}Let $(\ggo, J, g)$ be a K\"ahler Lie algebra and assume 
that  $\hh$ is a ideal satisfying $J\hh = \hh^{\perp}$ and $\hh \cap J\hh = 0$ 
(that is $\hh$ is $\omega$-lagrangian as in (\ref{m11})) then for 
$x, y \in \hh$ it holds
\begin{itemize}

\item $\nabla _x y \in J\hh$;

\item $\nabla_{Jx}Jy \in J\hh$;

\item $\nabla _x Jy \in \hh$; and $\nabla _{Jx} y \in \hh$

\end{itemize}
Thus the subgroup corresponding to $J\hh$ on the Lie group $G$ is totally geodesic.
\end{prop}

\begin{prop} \label{l2} Let $(\ggo, J, g)$ be a K\"ahler Lie algebra and 
assume that  $\hh$ is a abelian ideal satisfying $J\hh = \hh=\hh^{\perp}$. 
Thus $\ggo=\hh \ltimes \kk$ with $J\kk \subset \kk$. Then it holds: 
\begin{itemize}

\item $\nabla _z y \in \hh$ for all  $y \in \hh$, and $z \in \ggo$; in particular $\nabla_x y =0$ for all $x,y\n \hh$

\end{itemize}
Thus the normal subgroup $H$ corresponding to the ideal $\hh$ on the Lie group $G$ is totally geodesic.
\end{prop}

The proofs of the previous two propositions follow from the Koszul formula for
the Levi Civita connection and the features announced in Propositions
(\ref{m11}) and (\ref{m22}).

Recall that a pseudo  metric on a Lie algebra $\ggo$ is called a {\it Walker} metric if there exists a
null and parallel subspace $W\subset \ggo$, i.e. there is $W$ satisfying $g(W,W)=0$ and
$\nabla_y W\subset W$ for all $y$ (see \cite{Wa}). The previous proposition show examples
of Walker metrics in dimension four (compare with \cite{Mt}). In fact $W=\hh$ satisfies $\nabla_y W \subset W$ for all $y\in \ggo$.

\begin{cor} \label{cowa} The neutral metrics on the K\"ahler Lie algebras of Proposition 
(\ref{m22}) are Walker.
\end{cor}

The curvature tensor $R(x,y)$ and the Ricci tensor $ric(x,y)$ are respectively defined by:
$$R(x,y) = [\nabla_x, \nabla_y]-\nabla_{[x,y]} \qquad ric(x,y)=-\sum_i\varepsilon_i g(R(x,v_i)y,v_i)$$
where $\{v_i\}$ is a frame field on $\ggo$ and $\varepsilon_i$ equals $g(v_i,v_i)$. The left invariant property allows to speak in the following setting. We  say that the metric is flat if $R \equiv 0$, and similar we get Ricci flat or completeness of $\nabla$ at the level of $\ggo$.

It is clear that the existence of flat or non flat pseudo K\"ahler metrics is a property which is invariant under complex isomorphisms, i.e. if $J$ and $J'$ are equivalent complex structures then there exists a flat (resp. non flat) pseudo K\"ahler metric for $J$ if and only there exists such a metric for $J'$.

\begin{thm}\label{unif} Let $\ggo$ be a unimodular four dimensional K\"ahler Lie algebra with pseudo K\"ahler metric g, then $g$ is flat and its Levi Civita connection is complete.
\end{thm}
\begin{proof} Among the K\"ahler Lie algebras of (\ref{kahler}) the unimodular ones are $\RR \times \hh_3$ and $\RR\times\ee(2)$.

 In the first case, $\RR\times \hh_3$, any pseudo K\"ahler metric has the 
 form $g= a_{12}( e^1 . e^1+ e^2 . e^2) + a_{13} (e^2 . e^3- e^1.e^4) + 
 a_{14}( e^1 . e^3+ e^2 . e^4)$ and the corresponding  Levi Civita connection is
$$\nabla_z y= \frac1 {\varepsilon}[(\alpha y_1 +\beta y_2) e_3 +(\alpha y_2 -\beta y_1) e_4]$$
where $\varepsilon= a_{13}^2+a_{14}^2$, $\alpha= -a_{13}(a_{13}z_2+a_{14} z_1)$ and $\beta= a_{14}(a_{14}z_1+a_{13}z_2)$.

For the Lie algebra $\RR\times\ee(2)$, any pseudo K\"ahler metric is $g= a_{14}( e^1 . e^1- e^4 . e^4) + a_{23}( e^2 . e ^2+ e^3 . e^3)$ and the corresponding Levi Civita connection is
$$\nabla_z y= z_1 y_3 e_2 - z_1y_2 e_3$$
In both cases the connection $\nabla$ is complete: for $\RR \times \hh_3$ the geodesic equations follows:
$$ x_1'=0,\, x_2'=0,\, x_3'=\frac 1{\varepsilon} (\alpha x_1+\beta x_2),\, x_4'= \frac 1{\varepsilon}( \alpha x_2 -\beta x_1)$$
and for $\RR \times \ee(2)$:
$$ x_1'=0,\, x_2'= x_1 x_3,\, x_3'=-x_1 x_2,\, x_4'=0$$
whose solution for a given initial condition are defined in $\RR$. In both cases $\nabla_{[x,y]} \equiv 0$ and since $\nabla_x \nabla_y = \nabla_y \nabla_x$, the curvature tensor  vanishes which implies that $g$ is flat.
\end{proof}

In the non definite case Ricci flat metrics do not need to be flat. Known counterexamples for neutral metrics are provided by hypersymplectic Lie algebras.


In the left invariant case, Lie algebras admitting hypersymplectic structures
are examples of K\"ahler Lie algebras, with some extra structure. 
In fact, if $(\ggo,J,E,g)$ is a hypersymplectic Lie algebra, then $J$ is a complex structure, $E$ a product structure  that anticommutes with $J$, and $g$ is a compatible metric  such that the associated 2-forms are closed. Furthermore $\ggo$ admits a splitting as vector subspaces  $\ggo=\ggo_+ \oplus \ggo_-$ of subalgebras of $\ggo$, with $J\ggo_+=\ggo_-$. Then $\ggo_+$ carries a flat torsion free connection $\nabla^+$ compatible with a symplectic form $\omega_+$, and similarly, $\ggo_-$ carries a flat torsion free connection $\nabla_-$ and a compatible symplectic form $\omega_-$. Both symplectic forms are related by $\omega_+(x,y)=\omega_-(Jx,Jy)$ for  $x,y \in \ggo_+$ (see for instance \cite{Ad}).

Such metric is neutral and Ricci flat. We find  more examples of Ricci flat 
metrics than the hypersymplectic ones  in the four dimensional case. 


 Hypersymplectic four dimensional Lie algebras were classified in \cite{Ad}.  
 Aside from the abelian Lie algebra there are only three Lie algebras which 
 admit a hypersymplectic structure: $\RR \times \hh_3$, $\rr_{4,-1,-1}$ and 
 $\dd_{4,2}$. In Theorem (\ref{unif}) it was proved that the Lie algebra $\RR \times \ee(2)$ is 
 flat and it does not admit hypersymplectic structures \cite{Ad}. In the 
 following theorem we shall complete the list of K\"ahler Lie algebras 
 $(\ggo, J, g)$ whose pseudo K\"ahler metric is Ricci flat.

 \begin{remark} It is known that for a given complex product structure on a four dimensional Lie algebra there is only one compatible metric, up to a non zero constant (see for instance \cite{Ad}).
 \end{remark}

\begin{thm} \label{t2}Let  $(\ggo,J)$ be  a non unimodular four dimensional K\"ahler Lie algebra with pseudo K\"ahler metric $g$ which is Ricci flat.  Then  $(\ggo, J)$ is isomorphic either to  $(\rr_{4,-1,-1}, J)$, $(\dd_{4,2}, J_2)$, $(\aff(\CC),  J_2)$ . Moreover these Lie algebras have flat metrics and also  Ricci flat but non flat metrics.
\end{thm}
\begin{proof} For each one of these Lie algebras we will exhibit all pseudo 
K\"ahler metrics in its matricial representation, and the computations prove 
that they are  Ricci flat. In particular for all $s$ such that $s=0$ the 
corresponding metric is flat.
$$
\begin{array}{ll}
\begin{array}{l}
\aff(\CC):\\
J_2e_2=e_1,\, J_2e_3=e_4 \\ \\
\left( \begin{matrix}
-s & 0 & a_{14} & -a_{13}  \\
0 & -s & -a_{13} & -a_{14} \\
 a_{14} & -a_{13} & 0 & 0 \\
-a_{13} & -a_{14} & 0 & 0
\end{matrix}
\right)
 \end{array}
& {\begin{array}{l}
 \nabla_Z Y  =   (-z_1 y_1 + z_2 y_2) e_1 - (z_2 y_1 +  z_1 y_2) e_2 \\
\quad \quad + (\frac{s}{\varepsilon}\alpha y_1+\frac{s}{\varepsilon}\beta y_2 + z_1 y_3 - z_2 y_4) e_3 + \\
\quad \quad +  (\frac{s}{\varepsilon}\beta  y_1 -\frac{s}{\varepsilon}\alpha y_2 + z_2 y_3 + z_1 y_4) e_4\\
\varepsilon  =  a_{13}^2+a_{14}^2 \\
\alpha  =  -a_{14}z_1 + a_{13}z_2 \\
\beta  =  a_{13} z_1 + a_{14}z_2 \\
R(X,Y)Z =  2\frac{(x_1y_2-x_2y_1)}{\varepsilon}[(a_{13}z_1+ \\ \qquad \qquad \quad a_{14}z_2) e_3 + (a_{14}z_1-a_{13}z_2)e_4] \\
ric(X,Y) =0\\
g(R(v,w)w,v)=-s(v_1w_2-v_2w_1)^2
\end{array}} \\ \\ 
\begin{array}{l}
\rr_{4,-1,-1}:\\
Je_4=e_1,\, J e_2=e_3 \\ \\ 
\left( \begin{matrix}
-s & a_{13} & -a_{12} & 0 \\
a_{13} & 0 & 0 & -a_{12} \\
-a_{12} & 0 & 0 & -a_{13} \\
0 & -a_{12} & -a_{13} & -s
\end{matrix}
\right)
\end{array}
& \begin{array}{l}
\nabla_Z Y =\frac1{\varepsilon}[\varepsilon z_4 y_1 e_1 + (s \alpha y_1 -
\varepsilon z_4 y_2 +s \beta y_4) e_2\\
 +(s \beta y_1-\varepsilon z_4 y_3 -s \alpha y_4) e_3 + \varepsilon z_4 y_4 e_4]\\
 \varepsilon= a_{12}^2+a_{13}^2\\
\alpha = a_{12} z_1 + a_{13} z_4 \\
\beta = a_{13} z_1 - a_{12} z_4 \\
R(X,Y)Z = \frac{3s(x_1y_4-x_4y_1)}{\varepsilon}[(a_{12} z_1+a_{13}z_4)e_2 \\
\qquad \qquad \quad + (a_{13} z_1-a_{12}z_4)e_3]\\
ric(X,Y) = 0\\
g(R(v,w)w,v) = s(v_4w_1-v_1w_4)^2
\end{array} \\ \\ \end{array}$$
$$\begin{array}{ll}
\begin{array}{l}
\dd_{4,2}:  \\J_1 e_2= e_4,\, J_1 e_1=e_3 \\ \\
\left( \begin{matrix}
0 & a_{14}  & 0 & 0 \\
 a_{14} & s  & 0 & 0\\
0 & 0 & 0 & a_{14} \\
0 & 0 & a_{14} & s
\end{matrix}
\right)
\end{array}
& \begin{array}{l}
\nabla_Z Y  =  (z_4y_1+\frac{s}{a_{14}}z_4y_2+(-z_1+\frac{s}{a_{14}}z_2)y_4)e_1- \\
- z_4y_2e_2+ ((z_1-\frac{s}{a_{14}}z_2)y_2+ z_4y_3+\frac{s}{a_{14}}z_4y_4)e_3 -z_4 e_4\\
R(X,Y)Z = 3 \frac{s}{a_{14}}(x_4 y_2-x_2 y_4 ) [z_4e_1-z_2 e_3] \\
ric(X,Y) = 0\\
g(R(v,w)w,v)=-3(v_4w_2-w_4v_2)^2
\end{array} \\
\end{array} $$

The other K\"ahler Lie algebras do not admit Ricci flat metrics (see  results of Proposition (\ref{t3}) and Theorem (\ref{t4})).

Notice that in all cases the commutator is a totally geodesic submanifold.
Moreover in $\aff(\CC)$ we have $\nabla_{\ggo'}\ggo' =0$, and in the other 
cases $\nabla_{\ggo'}\ggo'\subset span\{e_3\}$ for any $s$. If $s = 0$ then in 
$\rr_{4,-1,-1}$ we get that the Levi Civita connection restricted to the 
commutator is always zero. Furthermore the abelian ideal $\hh$ is flat  
where $\hh=\ggo'$ in $\aff(\CC)$, $\hh=span\{e_2,e_3\}$ in $\rr_{4,-1,-1}$ and 
$\hh=span\{e_1, e_3\}$ in $\dd_{4,2}$ (see (\ref{cowa}).
\end{proof}

\begin{remark} Among these Ricci flat metrics there are examples of complete and non complete metrics \cite{Ad}.
\end{remark}

\medspace

An {\it Einstein metric} $g$ is proportional to its corresponding Ricci tensor, i.e. $g(x,y) = \nu \, ric(x,y)$ for all $x, y \in \ggo$ and $\nu$ be a real constant. We shall determine Einstein K\"ahler metrics in the four dimensional case.

\begin{prop} \label{t3} Let $(\ggo, J, g)$ be a K\"ahler Lie algebra with Einstein metric 
$g$. Then if $g$ is non Ricci flat, g is a pseudo K\"ahler metric corresponding to one of the following Lie algebras:
$$\begin{array}{lll}
\aff(\RR) \times \aff(\RR) & J & g= \alpha( e^1 . e^1 + e^2 . e^2 + e^3 . e^3 + e^4 . e^4)\\
\aff(\CC) &  J_1 & g= \alpha( e^1 . e^1 - e^2 . e^2 + e^3 . e^3 - e^4 . e^4)\\
\dd_{4,1/2} & {\begin{array}{l}
J_1,\\
 J_2\end{array}} & {\begin{array}{l} \vspace{-.27cm}\\
 g= \alpha( e^1 . e^1 + e^2 . e^2 +  e^3 . e^3 + e^4 . e^4)\\
 g= \alpha( e^1 . e^1 + e^2 . e^2 -  e^3 . e^3 - e^4 . e^4)
 \\ \vspace{-.35cm}\end{array}} \\
\dd_{4,\delta}' & {\begin{array}{l}
J_1, J_3 \\
J_2, J_4
\end{array}}
& {\begin{array}{l}
\vspace{-.35cm}\\
g= \alpha( e^1 . e^1 + e^2 . e^2 +  \delta ( e^3 . e^3 + e^4 . e^4))\\
g= \alpha( e^1 . e^1 + e^2 . e^2 -  \delta ( e^3 . e^3 + e^4 . e^4))\\
\vspace{-.27cm}
\end{array}}
\end{array}
$$
In all cases $\alpha\neq 0$.
\end{prop}
\begin{proof} We shall exhibit  all pseudo K\"ahler metrics on these Lie algebras. We compute the corresponding  Levi Civita connection, curvature and Ricci tensor on each case.
$$
\begin{array}{ll}
\begin{array}{l}
\aff(\RR) \times \aff(\RR):\\
Je_1 = e_2, \, Je_3 = e_4 \\ \\
 \left( \begin{matrix}
a_{12} & 0 & 0 & 0 \\
0 & a_{12} & 0 & 0 \\
 0 & 0 & a_{34} & 0 \\
0 & 0 & 0 & a_{34}
\end{matrix}
\right)
\end{array}
& \begin{array}{l}
\nabla_Z Y = z_2(y_2 e_1 - y_1 e_2)+z_4(y_4 e_3 - y_3 e_4)\\
R(X,Y) = - \nabla_{[X,Y]}\\
ric(X,Y) = -x_1y_1 - x_2 y_2 - x_3y_3 - x_4 y_4 \\
g(R(v,w)w,v) = -a_{12}(v_2w_1-v_1w_2)^2 -a_{34}(v_3w_4-v_4w_3)^2
 \end{array} \end{array}
$$
Clearly if $a_{12}=a_{34} \neq 0$, then the corresponding metric is Einstein.
$$
\begin{array}{ll}
\begin{array}{l}
\aff(\CC):\\
J_1 e_1 = e_3,\, J_1 e_2 = e_4 \\ \\
\left( \begin{matrix}
a_{13} & a_{14} & 0 & 0 \\
a_{14} & -a_{13} & 0 & 0 \\
0 & 0 & a_{13} & a_{14} \\
0 & 0 & a_{14} & -a_{13}
\end{matrix}
\right)
\end{array}
& \begin{array}{l} 
\nabla_Z Y   =  (z_3 y_3 - z_4 y_4) e_1  + (z_4 y_3 + z_3 y_4) e_2-\\
\quad - (z_3 y_1 - z_4 y_2) e_3 -(z_4 y_1 + z_3 y_2) e_4 \\
R(X,Y)=-\nabla_{[X,Y]}\\
ric(X,Y) = 2(-x_1 y_1+x_2 y_2-x_3 y_3+x_4 y_4)\\
g(R(v,w)w,v)= -a_{13}(\alpha^2-\beta^2) - 2 a_{14}\alpha \beta\\
\alpha= v_1w_3-v_3w_1+v_4w_2-v_2w_4\\
\beta= w_4v_1-v_4w_1+v_2w_3-v_3w_2
 \end{array} \end{array}
$$
Therefore when $a_{14}=0$ and $a_{13}\neq 0$ the corresponding metric is Einstein. 

$$\begin{array}{ll}
\begin{array}{l}
\dd_{4,\frac12}: \\
J_ 1e_1 = e_ 2, \, J_1 e_4 = e_3 \\
\\ a_{12}
\left(\begin{matrix}
1 & 0 & 0 & 0 \\
0 & 1 & 0 & 0 \\
0 & 0 & 1& 0 \\
0 & 0 & 0 & 1
\end{matrix}
\right)
\end{array}
& \begin{array}{l}
\nabla_Z Y  =  \frac12(z_3y_2+z_2y_3-z_1y_4)e_1+ \\ \qquad \quad + \frac12(-z_3y_1-z_1y_3-z_2y_4)e_2+\\
\qquad \quad + {[ \frac12(-z_2y_1+z_1y_2)-z_3y_4 ]e_3} + \\
\qquad \quad + { [\frac12(z_1y_1+z_2y_2)+z_3y_3]e_4} \\
R(X,Y)Z  =  [(\alpha -\frac12 \eta) z_2- \frac14(\nu +\gamma)z_3 + \frac14(\theta + \beta)z_4] e_1 \\
\qquad \quad  +[(\frac12 \eta-\alpha)z_1+\frac14(\theta+\beta)z_3 +(-\frac14\nu+\frac12 \gamma)z_4]e_2\\
\qquad \quad + [\frac14 (-\nu +\gamma)z_1-\frac14(\theta+\beta)z_2+(\eta -\frac12 \alpha)z_4]e_3\\
\qquad  \quad +[-\frac14(\theta +\beta)z_1+\frac14(\nu-\gamma)z_2+ (\frac12\alpha -\eta)z_3]e_4 \\
\alpha= x_2y_1 -x_1y_2, \beta  = x_4y_3-x_3y_4, \gamma = x_4y_2-x_2y_4 \\
 \eta = x_4y_3-x_3y_4, \nu = x_3y_1-x_1y_3, \theta = x_3y-2-x_2y_3\\
ric = -\frac32 g\\
g(R(v,w)w,v)=a_{12}(\alpha^2+\eta^2-\eta\alpha-\frac14(\beta +\theta)^2+\frac14(\nu-\gamma)^2
\end{array} \\
\begin{array}{l}
\mbox{ for } J_2: \\
J_2 e_1=e_2, \, J_2 e_4 = e_3 \\
\\
a_{12} \left(\begin{matrix}
-1 & 0 & 0 & 0 \\
0 & -1 & 0 & 0 \\
0 & 0 & 1& 0 \\
0& 0 & 0 & 1
\end{matrix}
\right) 
\end{array}
& \begin{array}{l} 
\nabla_Z Y  =  -\frac12(z_3y_2+z_2y_3+z_1y_4)e_1+ \\ \qquad \quad + \frac12(z_3y_1+z_1y_3-z_2y_4)e_2+\\
 \qquad \quad + { [\frac12(-z_2y_1+z_1y_2)-z_3y_4]} e_3 + \\
 \qquad \quad + {[-\frac12(z_1y_1+z_2y_2)+z_3y_3]} e_4 \\
R(X,Y)Z  = \\
\qquad  {[(-\alpha +\frac12 \eta) z_2+\frac14(\nu +\gamma)z_3 + \frac14(-\theta + \beta)z_4] e_1} \\
\qquad { +[(-\frac12 \eta+\alpha)z_1+\frac14(\theta-\beta)z_3 +\frac14 (\nu +\gamma)z_4]e_2}\\
\qquad +{ [\frac14 (\nu +\gamma)z_1+\frac14(\theta-\beta)z_2+(\eta -\frac12 \alpha)z_4]e_3}\\
\qquad +{[\frac14(-\theta +\beta)z_1+\frac14(\nu+\gamma)z_2+ (\frac12\alpha -\eta)z_3]e_4} \\
ric = - \frac32 g\\
g(R(v,w)w,v)=a_{12}[(\alpha^2+\eta^2-\eta \alpha)-\frac14(\nu + \gamma)^2-\frac14(\beta-\theta)^2]
\end{array} \end{array}
$$ 
where $\alpha, \beta, \gamma,  \eta, \nu, \theta$ are as above. 
Therefore any pseudo  K\"ahler metric is Einstein. 

On $\dd'_{4,\delta}$
 we have four non equivalent complex structures compatible with the same symplectic structure $\omega = a_{12} (e^1\wedge e^2-\delta e^3\wedge e^4) $, with $a_{12} \ne 0$:
$$  J_1 e_1 = e_2 \qquad J_1e_4 = e_3 \qquad J_2 e_1 = e_2 \qquad J_2e_4 = -e_3 $$
$$  J_3 e_1 = -e_2 \qquad J_3e_4 = -e_3 \qquad J_4 e_1 = -e_2 \qquad J_4e_4 = e_3 $$

 The corresponding  pseudo-Riemannian K\"ahler metrics are:

$$\text{ for $J_1$: } a_{12}
\left(\begin{matrix}
1 & 0 & 0 & 0 \\
0 & 1 & 0 & 0 \\
0 & 0 & \delta& 0 \\
0& 0 & 0 & \delta
\end{matrix}
\right)\qquad \text{ for $J_2$: }
a_{12} \left(\begin{matrix}
1 & 0 & 0 & 0 \\
0 & 1 & 0 & 0 \\
0 & 0 & -\delta & 0 \\
0& 0 & 0 & -\delta
\end{matrix}
\right)
$$

$$\text{ for $J_3$: } -a_{12}
\left(\begin{matrix}
1 & 0 & 0 & 0 \\
0 & 1 & 0 & 0 \\
0 & 0 & \delta& 0 \\
0& 0 & 0 & \delta
\end{matrix}
\right)\qquad \text{ for $J_4$: }
-a_{12} \left(\begin{matrix}
1 & 0 & 0 & 0 \\
0 & 1 & 0 & 0 \\
0 & 0 & -\delta & 0 \\
0& 0 & 0 & -\delta
\end{matrix}
\right)
$$

We investigate two cases. The Levi-Civita connection for $g_1$ is:

$$\begin{array}{rcl}
\nabla_Z Y & = & [(z_4+\frac{\delta}2 z_3)y_2+\frac{\delta}2 (z_2y_3-z_1y_4)]e_1 \\
& & +[(\frac{\delta}2z_3+z_4)y_1-\frac{\delta}2(z_1y_3+z_2y_4)]e_2\\
& & +[\frac12(-z_2y_1+z_1y_2)-\delta z_3y_4]e_3 \\
& & + [\frac12(z_1y_1+z_2y_2)+\delta z_3y_3]e_4
\end{array}
$$
The curvature tensor is
$$\begin{array}{rcl}
R(X,Y)Z & = & [\delta(\alpha -\frac{\delta}2 \eta) z_2+\frac{\delta^2}4((\nu -\gamma)z_3 + (\theta + \beta)z_4] e_1 \\
& & +[\delta(\frac{\delta}2 \eta-\alpha)z_1+\frac{\delta^2}4((\theta+\beta)z_3 +(-\nu +\gamma)z_4)]e_2\\
& & + [\frac{\delta}4 ((-\nu +\gamma)z_1-(\theta+\beta)z_2)+\delta(\delta \eta -\frac12 \alpha)z_4]e_3\\
&& +[\frac{\delta}4(-(\theta +\beta)z_1+(\nu-\gamma)z_2)+ \delta(\frac12\alpha -\delta\eta)z_3]e_4
\end{array}
$$

$$
\alpha= x_2y_1 -x_1y_2 \quad \beta  = x_4y_3-x_3y_4 \quad \gamma = x_4y_2-x_2y_4$$
$$ \eta = x_4y_3-x_3y_4 \quad \nu = x_3y_1-x_1y_3 \quad \theta = x_3y-2-x_2y_3$$

the Ricci tensor is
$$ric = -\frac32 \delta g_1$$
$$g(R(v,w)w,v)= -a_{12}(\delta(\alpha -\frac12 \delta\eta)\alpha +
\delta^2(\delta\eta-\frac12\alpha)\eta-\frac14\delta^2[(\gamma-\nu)^2+(
\theta+\beta)^2])$$

\

The Levi-Civita connection for $g_2$ is:

$$\begin{array}{rcl}
\nabla_Z Y & = & [(z_4-\frac{\delta}2 z_3)y_2-\frac{\delta}2 (z_2y_3+z_1y_4)]e_1 \\
& & +[(\frac{\delta}2z_3-z_4)y_1-\frac{\delta}2(z_1y_3-z_2y_4)]e_2\\
& & +[\frac12(-z_2y_1+z_1y_2)-\delta z_3y_4]e_3 \\
& & + [-\frac12(z_1y_1+z_2y_2)+\delta z_3y_3]e_4
\end{array}
$$
The curvature tensor is
$$\begin{array}{rcl}
R(X,Y)Z & = & [\delta(-\alpha +\frac{\delta}2 \eta) z_2+\frac{\delta^2}4((\nu +\gamma)z_3 + (-\theta + \beta)z_4] e_1 \\
& & +[\delta(-\frac{\delta}2 \eta+\alpha)z_1+\frac{\delta^2}4((\theta-\beta)z_3 +(\nu +\gamma)z_4)]e_2\\
& & + [\frac{\delta}4 ((\nu +\gamma)z_1-(\theta-\beta)z_2)+\delta(\delta \eta -\frac12 \alpha)z_4]e_3\\
&& +[\frac{\delta}4((-\theta +\beta)z_1+(\nu+\gamma)z_2)+ \delta(\frac12\alpha -\delta\eta)z_3]e_4
\end{array}
$$

$$
\alpha= x_2y_1 -x_1y_2 \quad \beta  = x_4y_3-x_3y_4 \quad \gamma = x_4y_2-x_2y_4$$
$$ \eta = x_4y_3-x_3y_4 \quad \nu = x_3y_1-x_1y_3 \quad \theta = x_3y-2-x_2y_3$$

the Ricci tensor is
$$ric = \frac32 \delta g$$
$$g(R(v,w)w,v)=-a_{12}[(-\alpha \delta+\frac12\delta^2\eta)\alpha+\delta^2(-\frac12\alpha +\delta \eta)\eta+\frac14\delta^2((\nu +\gamma)^2+(\beta -\theta)^2)]
$$

The proof will be completed with the results of the Theorem (\ref{t4}), 
proving that there is no more Einstein metrics.
\end{proof}

 We are in conditions to finish this geometric study with the characterization of four
 dimensional K\"ahler Lie algebras which are not Einstein.
 
 \begin{thm} \label{t4} Let $(\ggo, J, g)$ be a K\"ahler Lie algebra. If $\ggo$  does not
 admit an Einstein K\"ahler metric then $\ggo$ is isomorphic to $\RR^2\times
 \aff(\RR)$, $\rr'_{4,0,\delta}$, $\dd_{4,1}$.
 \end{thm} 
 \begin{proof} The previous propositions show all examples of Lie algebras
 admitting Einstein K\"ahler pseudo metrics. We shall show that the Lie algebras
 $\RR^2\times\aff(\RR)$, $\rr'_{4,0,\delta}$, $\dd_{4,1}$ do not admit Einstein
 K\"ahler metrics in a case by case study. 
$$\begin{array}{ll}
\begin{array}{l}
\RR^2 \times \aff(\RR):\, 
J e_1= e_2,\, J e_3 = e_4\\
\nabla_Z Y = z_2(y_2 e_1 - y_1 e_2) \\
R(X,Y) = - \nabla_{[X,Y]} \\
ric(X,Y) = -x_1y_1 - x_2 y_2 \\
g(R(v,w)w,v) = -a_{12}(v_2w_1-v_1w_2)^2
\end{array}
& \left( \begin{matrix}
a_{12} & 0 & 0 & 0 \\
0 & a_{12} & 0 & 0 \\
 0 & 0 & a_{34} & 0 \\
0 & 0 & 0 & a_{34}
\end{matrix}
\right) \\ \\ 
\begin{array}{l}
\rr'_{4,0,\delta}: \text{ for } J_1: J_1 e_4= e_1 \, J e_2= e_3\\
\nabla_Z Y =- z_1 y_4 e_1 +  \delta  z_4 y_3 e_2 - z_4 y_2 e_3 +  z_1 y_1 e_4\\
R(X,Y) = -\nabla_{[X,Y]} \\
ric(X,Y) = -x_1 y_1 - x_4 y_4\\
g(R(v,w)w,v)=a_{14}(v_4w_1-w_4v_1)^2
\end{array}
 & \left( \begin{matrix}
-a_{14} & 0 & 0 & 0 \\
 0 & a_{23}  & 0 & 0 \\
0 & 0 & a_{23} & 0 \\
0 & 0  & 0 & -a_{14}
\end{matrix}
\right) \\ \\ \begin{array}{l}
\text{ for } J_2:\, Je_4=e_1\, Je_2=-e_3 \\
\nabla_Z Y =- z_1 y_4 e_1 -  \delta  z_4 y_3 e_2 + z_4 y_2 e_3 +  z_1 y_1 e_4\\
R(X,Y) = -\nabla_{[X,Y]} \\
ric(X,Y) = -x_1 y_1 - x_4 y_4\\
g(R(v,w)w,v)=a_{14}(v_4w_1-w_4v_1)^2
\end{array}
& \left( \begin{matrix}
a_{14} & 0 & 0 & 0 \\
 0 & a_{23}  & 0 & 0 \\
0 & 0 & a_{23} & 0 \\
0 & 0  & 0 & a_{14}
\end{matrix}
\right)\end{array}$$
$$\begin{array}{ll} 
\begin{array}{l}
\dd_{4,1}: 
\nabla_Z Y =- z_1 y_4 e_1 -(z_3 y_1+z_1y_3) e_2 +\\
\qquad \quad  \qquad + (z_1 y_2-z_3y_4) e_3 +  z_1 y_1 e_4\\
R(X,Y) = -\nabla_{[X,Y]}\\
ric(X,Y) = -2(x_1 y_1 + x_4 y_4)\\
g(R(v,w)w,v) = -\alpha(a_{14}\alpha-2 \beta a_{12})\\
\alpha=v_4w_1-v_1w_4\\
\beta=v_1w_2-w_1v_2+v_4w_3-w_4v_3
\end{array}
& \begin{array}{l}
J e_1= e_4,\, J e_2=e_3 \\
\left( \begin{matrix}
a_{14} & 0 & -a_{12} & 0 \\
 0 & 0  & 0 & a_{12} \\
-a_{12} & 0 & 0 & 0 \\
0 & a_{12}  & 0 & a_{14}
\end{matrix}
\right) \end{array}\\ \\ \end{array}$$ 
Finally notice that the Lie algebra $\dd_{4,2}$ admits two non equivalent
complex structures, one of them admits a compatible Einstein pseudo metric. But
for the other one $J_2$ this is not the case as the following computations show.
$$\begin{array}{ll}
\begin{array}{l}
\nabla_Z Y  = (\frac{a_{23}}{a_{14}}(z_3 y_2 +z_2y_3) -2z_1y_4)e_1+ \\
\qquad +(\frac12(-z_3 y_1-z_1y_3)+z_2y_4) e_2+ \\
\qquad +(\frac12(-z_2y_1+z_1y_2)-z_3y_4)e_3+ \\
\qquad +(\frac12z_1 y_1-\frac{a_{23}}{2a_{14}}(-z_2y_2+z_3y_3) e_4 \\
R(X,Y)Z = -\frac{a_{23}}{a_{14}}{[(\alpha + \eta) z_2+(\frac12\nu +\gamma)z_3 + (2\theta + }\\  {\qquad 4\beta)z_4]} e_1
+{[(\frac12 \eta-\alpha)z_1+\frac{a_{23}}{a_{14}}(-\theta+\beta)z_3 +\frac12 \nu +\gamma)z_4]}e_2\\
\qquad + {[(-\frac14 \nu -\frac12 \gamma)z_1+(\frac{a_{23}}{a_{14}}\theta-\beta)z_2+(\eta -\frac12 \alpha)z_4]}e_3\\
\qquad + {\frac{a_{23}}{a_{14}}[(-\frac12 \theta -\beta)z_1+(-\theta -\frac12 \eta)z_2+ \frac12(\alpha -\eta)z_3]}e_4 \\
\alpha= x_2y_1 -x_1y_2, \beta  = x_4y_3-x_3y_4, \gamma = x_4y_2-x_2y_4\\
\eta = x_4y_3-x_3y_4, \nu = x_3y_1-x_1y_3, \theta = x_3y-2-x_2y_3\\
ric(X,Y)=-6x_4y_4 -\frac32 x_1y_1
\end{array}
& \begin{array}{l}
J_2 e_4 = -2e_1,\, J_2 e_2=e_3 \\ \\
\left( \begin{matrix}
\frac12 a_{14} & 0 & 0  & 0 \\
 0 & a_{23}  & 0 & 0 \\
0 & 0 &  a_{23} & 0\\
0 & 0 & 0 & 2a_{14}
\end{matrix}
\right) \end{array} \\ \end{array}
$$
\end{proof}

\begin{cor} Let $(\ggo,J)$ be a four dimensional K\"ahler Lie algebra. Then
 the commutator is  totally geodesic.
\end{cor}
 \begin{proof} It follows from the Levi Civita connection computed at the corresponding elements in the commutator.
 \end{proof}

\section{A picture in global coordinates}

In this section we shall write the pseudo K\"ahler metrics in global complex 
coordinates (the real expression can also be done with the information we
present in the following paragraphs). The following table summarizes the results. In the first column we
write the corresponding Lie algebra, the invariant complex structure and the homogeneous
complex manifold according to \cite{Sn} and \cite{O1}. In the second column we
present left invariant 1-forms and the metric in terms of complex coordinates.

$$\begin{array}{ll}
{\begin{array}{c} \vspace{-.27cm}\\
\RR \times \hh_3 \\ Jv_1=v_2, \, Jv_3=v_4 \\ \CC^2 \\ \vspace{-.27cm}\end{array}} & {\begin{array}{c} 
v^1= dx, \, v^2= dy, \, v^3=dz +\frac y2 dx -\frac x2 dy,\, v^4= dt\\
\mbox{ with } u=v_1 + iv_2, \, w= v_3 + i v_4\\
g=a_{12} du d\overline{u}+(a_{14}-ia_{13}) du d\overline{w}+(a_{14}+ia_{13})
d\overline{u}dw \\ \mbox{ flat (\ref{unif}) }\end{array}} \\
{\begin{array}{c} \vspace{-.27cm}\\
\RR^2 \times \aff(\RR) \\ Jv_1=v_2, \, Jv_3=v_4 \\ \CC\times \HH \\ 
\vspace{-.27cm}\end{array}} &
{\begin{array}{c} v^1= dt, \, v^2= e^{-t}dx , \, v^3=dy,\, v^4= dz \\ 
\mbox{ with } u=v_1 + iv_2, \, w= v_3 + i v_4 \\
g=a_{12} du d\overline{u}+a_{34}dw d\overline{w} \end{array}} \\
{\begin{array}{c} \vspace{-.27cm}\\
\RR \times \ee(2) \\ Jv_1=v_4,\,Jv_2=v_3 \\ \CC^2 \\ 
\vspace{-.27cm}\end{array}} &
{\begin{array}{c} v^1= dt, \, v^2= \cos{t}dx + \sin{t} dy , \, v^3=\sin{t} dx + \cos{t} dy ,\, 
v^4= dz \\ 
\mbox{ with } u=v_1 + iv_4, \, w= v_2 + i v_4 \\
g=a_{14} du d\overline{u}+a_{23}dw d\overline{w} \end{array}}\\
{\begin{array}{c} \vspace{-.27cm}\\
\aff(\RR)\times \aff(\RR) \\ Jv_1=v_2,\,Jv_3=v_4 \\ \HH \times \HH \\ 
\vspace{-.27cm}\end{array}} &
{\begin{array}{c} v^1= dx, \, v^2= e^{-x}dy, \, v^3=dz,\, v^4= e^{-z} dt\\ 
\mbox{ with } u=v_1 + iv_2, \, w= v_3 + i v_4 \\
g=a_{12} du d\overline{u}+a_{34}dw d\overline{w} \\
\mbox{ Einstein if } a_{12}=a_{34} \neq 0 \, (\ref{t3}) \end{array}}\\
{\begin{array}{c} \vspace{-.27cm}\\
\aff(\CC) \\ J_1v_1=v_3,\,J_1v_2=v_4 \\ \CC^2 \\ 
\end{array}} &
{\begin{array}{c} v^1= dt, \, v^2= dz, \, v^3=e^{-t}(\cos z \, dx+ \sin z\, dy),\, 
 \\ 
v^4= e^{-t} (-\sin z \, dx + \cos z \, dy),\, \mbox{ with } u=v_1 + iv_3, \, w= v_2 + i v_4 \\
g_1=a_{13} (du^2+ d\overline{u}^2 + dw^2+ d\overline{w}^2)+a_{14}i(du^2-
d\overline{u}^2 + dw^2- d\overline{w}^2)\\
\mbox{ Einstein if } a_{14} = 0 \, (\ref{t3}) \end{array}}\\
{\begin{array}{c}  J_2v_1=-v_2,\, J_2 v_3=v_4\\
\CC^2 \\ \vspace{-.27cm} \end{array}} & {\begin{array}{c}  
\mbox{ with } u=v_1 + iv_2, \, w= v_3 + i v_4 \\
g_2=s du d\overline{u} + a_{14}(du dw + d\overline{u} d\overline{w})-i
a_{13} (du dw - d\overline{u}d\overline{w})\\
\mbox{ Ricci flat always and flat if } s = 0 \, (\ref{t2}) \end{array}}\\
{\begin{array}{c} \vspace{-.27cm}\\
\rr_{4,-1,-1} \\ Jv_4=v_1,\,Jv_2=v_3 \\ \CC\times \HH \\ 
\end{array} \vspace{-.27cm}} &
{\begin{array}{c} v^1= e^{-t} dx, \, v^2= e^{t}dy, \, v^3=e^{t}dz,\, v^4= dt\\ 
\mbox{ with } u=v_4 + iv_1, \, w= v_2 + i v_3 \\
g=-s du d\overline{u} -(a_{12}+i a_{13})du d\overline{w} -(a_{12}-i a_{13})
 d\overline{u} d{w}\\
\mbox{ Ricci flat always and flat if } s = 0 \, (\ref{t2}) \end{array}}\\
{\begin{array}{c} \vspace{-.27cm}\\
\rr'_{4,0,\delta} \\ J_1v_4=v_1,\,J_1v_2=v_3 \\ \CC\times \HH \\ 
\end{array}} &
{\begin{array}{c} v^1= e^{-t} dx, \, v^2= (\cos t dy + \sin t dz), \, 
\, v^4= dt\\ 
v^3=( -\sin t dy + \cos t dz),\,\mbox{ with } u=v_4 + iv_1, \, w= v_2 + i v_3 \\
g_1=-a_{14}du d\overline{u} +a_{23}dw d\overline{w}
 \end{array}}\\ \end{array} $$
$$\begin{array}{ll}
 {\begin{array}{c} J_2v_4=v_1,\,J_2v_2=-v_3 \\ \CC\times \HH \\ 
\end{array} } &
 {\begin{array}{c} 
 \mbox{ with } u=v_4 + iv_1, \, w= v_2 + i v_3 \\
 g_2=a_{14}du d\overline{u} +a_{23}dw d\overline{w}\\ \vspace{-.27cm}
 \end{array}}\\
 {\begin{array}{c} \vspace{-.27cm}\\
\dd_{4,1} \\ Jv_1=v_4,\,Jv_2=v_3 \\ \CC\times \HH \\ 
\end{array}} &
{\begin{array}{c} v^1= e^{-t} dx, \, v^2=dy, \, v^3=e^{-t} dz - \frac{x}2 e^{-t}
dy,\, v^4= dt\\ 
\mbox{ with } u=v_1 + iv_4, \, w= v_2 + i v_3 \\
g_1=a_{14}du d\overline{u} -i a_{12}(du d\overline{w}-d\overline{u}dw)
 \end{array}}\\
  {\begin{array}{c} \vspace{-.27cm}\\
\dd_{4,1/2} \\ J_1v_1=v_2,\,J_1v_2=v_3 \\ \DD^2\\ 
\end{array}} &
{\begin{array}{c} v^1= e^{-t/2} dx, \, v^2=e^{-t/2}dy, \, v^3=e^{-t} dz - 
\frac{x}2 e^{-t} dy,\, v^4= dt\\ 
\mbox{ with } u=v_1 + iv_2, \, w= v_4 + i v_3 \\
g_1=a_{12}(du d\overline{u} + dw d\overline{w})\\ \mbox{ Einstein (\ref{t3})}
\end{array}}\\
 {\begin{array}{c} 
J_2v_1=-v_2,\,J_2v_2=v_3 \\ ({\DD^2}^c)^0\\ \vspace{-.27cm}
\end{array}} &
{\begin{array}{c} 
\mbox{ with } u=v_1 + iv_2, \, w= v_4 + i v_3 \\
g_2=a_{12}(-du d\overline{u} + dw d\overline{w})\\ \mbox{ Einstein (\ref{t3})}\\
\vspace{-.27cm}
 \end{array}}\\
  {\begin{array}{c} \vspace{-.27cm}\\
\dd_{4,2} \\ J_1v_2=v_4,\,J_1v_1=v_3 \\ \CC \times \HH \\ 
\end{array}} &
{\begin{array}{c} v^1= e^{-2t} dx, \, v^2=e^{t}dy, \, v^3=e^{-t} dz - 
\frac{x}2 e^{-t} dy,\, v^4= dt\\ 
\mbox{ with } u=v_2 + iv_4, \, w= v_1 + i v_3 \\
g_1=s du d\overline{u} + a_{14}( du d\overline{w}+d\overline{u}dw)\\ 
\mbox{ Ricci flat always and flat if $s=0$ (\ref{t2})}\end{array}}\\
{\begin{array}{c}  J_2v_1=1/2 v_4,\,J_2v_2=v_3 \\ \CC \times \HH \\ 
\end{array}} &
{\begin{array}{c} 
\mbox{ with } u=\sqrt{2}/2 v_1 + i\sqrt{2}v_4, \, w= v_2 + i v_3 \\
g_2= a_{14} du d\overline{u}+a_{23} dw d\overline{w}\\ 
\vspace{ -.27cm}\\ \end{array}}\\
 {\begin{array}{c} \vspace{-.27cm}\\
\dd'_{4,\delta} \\ J_1v_1=v_2,\,J_1v_4=v_3 \\ \DD^2\\ 
\end{array}} &
{\begin{array}{c} v^1= e^{-\delta t/2} (\cos t dx-\sin t dy), \, 
v^2=e^{-\delta t/2}(\sin t dx + \cos t dy),\\ \, v^3=e^{-t} dz + 
{x} e^{-t\delta/2} (\sin 2t dx - \cos 2t dy),\, v^4= dt\\ 
\mbox{ with } u=v_1 + iv_2, \, w= v_4 + i v_3, \, 
g_1=a_{12} (du d\overline{u} +  \delta dw d\overline{w})\\ 
\mbox{ Einstein (\ref{t3})}\end{array}}\\
{\begin{array}{c}  J_2v_1= v_2,\,J_2v_4=-v_3 \\ ({\DD^2}^c)^0 \\ 
\end{array}} &
{\begin{array}{c} 
\mbox{ with }  u=v_1 + iv_2, \, w= v_4 + i v_3 \\
g_2= a_{12} (du d\overline{u}- \delta  dw d\overline{w}), 
\quad \mbox{ Einstein (\ref{t3})}\\ 
\vspace{ -.27cm}\\ \end{array}}\\
\end{array}$$

\section{Some generalizations}

Notice that the constructions of K\"ahler Lie algebras given in Propositions (\ref{m11}) and (\ref{m22}) can be done in higher dimensions.

\subsection{K\"ahler structures on affine Lie algebras} We shall generalize the K\"ahler structures on four dimensional affine Lie algebras. Furthermore we
 get higher dimensional examples of Ricci flat metrics, generalizing a
 K\"ahler structure on $\aff(\CC)$.
 
Let $A$ be an associative Lie algebra. Then $\aff(A)$ is the Lie algebra 
$A\oplus A$ with Lie bracket given by:
$$[(a,b)(c,d)]=(ac - ca , ad - cb)$$
An almost complex structure on $\aff(A)$ is defined by $K(a,b)=(b,-a)$ which is 
integrable and parallel for the torsion free connection $\nabla_{(a,b)}(c,d) = (ac, ad)$. 

Affine Lie algebras play an important role in the characterization of the solvable Lie algebras admitting an abelian complex structure \cite{BD2}.

Assume that $A$ is commutative and that $e_i$ i= 1$ \hdots $ n is a basis of
$A$. Let $v_i=(e_i,0)$ and $w_i=(0,e_i)$ be a basis of $\aff(A)$. Consider the
dual basis $v^i$ $w^i$ of $\aff(A)^{\ast}$ and define a non degenerate two form by
$\omega = \sum v^i \wedge w^i$, that is
$\omega((x_1,y_1)(x_2,y_2))=\sum_i(x_1^iy_2^i-x_2^iy_1^i)$. Indeed $\omega$ is $K$ invariant. Furthermore it
is closed. Denote with $u^i$ i=1,$\hdots,$n, the coordinates of $u\in A$. 
Let $(x_1,y_1), (x_2,y_2), (x_3,y_3)$ be elements on $\aff(A)$ . Then
$$\begin{array}{rcl}
d\omega((x_1,y_1),(x_2,y_2),(x_3,y_3)) & = & \omega(\sum_i(0,x_1^i
y_2^i-x_2^iy_1^i),(x_3,y_3))+\\
& & +  \omega( \sum_i((0, x_2^iy_3^i-x_3^iy_2^i),(x_1,y_1))+\\
& & + \omega( \sum_i((0, x_3^iy_1^i-x_1^iy_3^i),(x_2,y_2)) \\
& & -\sum_i [(x_1^i
y_2^i-x_2^iy_1^i)
x_3^i+(x_2^iy_3^i-x_3^iy_2^i)x_1^i+\\
& & + (x_3^iy_1^i-x_1^iy_3^i)x_2^i]\\
& = & 0
\end{array}
$$
\begin{prop} The Lie algebras $\aff(A)$  carry a K\"ahler structure for any
commutative algebra $A$.\end{prop}

\begin{example} In dimension four we find many examples of affine Lie algebras.
The list consists of the Lie algebras $\RR \times \hh_3$, $\RR^2 \times
\aff(\RR)$, $\aff(\RR) \times \aff(\RR)$, $\dd_{4,1}$ and $\aff(\CC)$ (see
\cite{BD2} for details).
\end{example}  

This K\"ahler structure does not necessarly induces a Ricci flat metric. See for
example  $\RR^2 \times \aff(\RR)$. 

Assume now that $A$ is a commutative complex algebra and consider $J$ to be the 
almost complex structure on $\aff(A)$ given by $J(a,b) =(-ia,ib)$. Let $\nabla$ be the connection on $\aff(A)$ given by $$\nabla_{(a,b)}(c,d)=(-ac,ad).$$
Then since $A$ is commutative $\nabla$ is torsion free. Furthermore the connection is flat. Indeed
$$R((a,b),(c,d))=\nabla_{[(a,b),(c,d)]}=0$$
and $J$ is parallel, that is $\nabla J=0$. We shall prove that $\nabla$ is a metric connection.

Take coordinates $u^i$  on $A \oplus 0$ and $w^i$ on $0 \oplus A$. 
Let $g$ be the (pseudo) metric on $\aff(A)$ defined by:
$$g((a,b),(c,d)) = \sum_i (du^i dw^i+d\overline{u}^i d\overline{w}^i)= \sum_i
Re\,(ad+bc)^i$$
then $\nabla$ is the Levi Civita connection of $g$. It is easy to verify that
$\nabla_{(a,b)}$ is skew symmetric with  respect to $g$.  

\begin{prop} The Lie algebras $\aff(A)$ are endowed with a neutral Ricci flat 
K\"ahler metric for a commutative complex algebra $A$.
\end{prop}

For a curve $(a(t), b(t))$ on $\aff(A)$ the geodesic equation related to the previous pseudo K\"ahler metric:  
$ - \nabla_{(a,b)}(a,b)=(a', b')$ gives rise the following system
$$\left\{ \begin{array}{lcr}
a' & = &  a^2 \\
b' & = & -ab
\end{array} \right.
$$
with non trivial solutions $ a= ( \kappa_1 - t)^{-1}$,  $b= \kappa_2(t-\kappa_1)$ 
for $\kappa_1, \kappa_2$ constants, showing that the metric is not complete 
except when $a=0$ and $b= \kappa$ is also constant.

\vskip .3cm

A Walker metric $g$  on a Lie algebra $\ggo$ (in the sense of \cite{Wa}) is characterized by the existence of a
subspace $W\subset \ggo$ satisfying:
\begin{equation} \label{walker} 
g(W,W)=0 \quad \mbox{ and } \quad \nabla_y W \subset W\quad 
\mbox{ for all } y \in \ggo\end{equation}
where $\nabla$ denotes the Levi Civita connection for $g$.

Since 
$$g([x,y],z)=g(\nabla_x y,z) + g(\nabla_y z, x)=0\quad \mbox{ for all } x,y, z
\in \ggo$$
then $W\subset W^{\perp}$. Thus if the dimension of $W$ is a half of the
dimension of $W$ then $W$ must be a subalgebra.

An hypersymplectic metric on a Lie algebra $\ggo$ is an example of a Walker
metric (see section 4). The following result explains how to construct
hypersymmplectic metrics from Walker K\"ahler metrics. The proof follows from
the previous observation and features of hypersymplectic Lie algebras (see
\cite{Ad} for instance).

\begin{prop} Let $g$ be a Walker K\"ahler metric on a Lie algebra $\ggo$ for
which a subspace $W\subset \ggo$ satisfies conditions (\ref{walker}) and assume that $\ggo= W \oplus
JW$ (direct sum as vector subspaces). Then $g$ is  an hypersymplectic metric on $\ggo$.
\end{prop}

The condition $W\oplus JW$ is necessary as proved by $(\aff(\CC),J_2, g_2)$. In
fact $g_2$ is a Walker K\"ahler metric but it is not hypersymplectic. The condition
 for $g$ of being K\"ahler is necessary as we see in the following example.
 
\begin{example} Consider on $\aff(\CC)$ the complex structure given by
$J(a,b)=(ia, ib)$ and let $g$ be  the metric defined by
$$g((a,b), (c,d))= Re(a\overline{d} + b \overline{c})$$ 
Then $g$ is compatible with $J$ and the Levi Civita connection for $g$ is
$$\nabla_{(a,b)}(c,d)=( -\frac12(a\overline{c}+c \overline{a}), a \frac{(d +
\overline{d})}2 + c\frac{(\overline{b}-b)}2 )$$
The complex structure $J$ is not parallel (see (\ref{affc})), hence this metric is not pseudo K\"ahler.
However the metric is Walker. In fact consider $W=\{(0,b)\}$, $b \in \CC$, and prove that
conditions (\ref{walker}) are satisfied (compare with \cite{Mt}).
\end{example}

The author is very grateful to I. Dotti Miatello for her invaluable comments.

\end{document}